\documentclass[preprint,1p,times]{elsarticle}
\usepackage{amsmath, amssymb, amsthm}
\usepackage{graphicx}
\usepackage{color,soul}
\usepackage{float}
\usepackage{lineno}
\usepackage{pgf}
\usepackage{epstopdf}

\journal{ }
\begin{document}

\begin{frontmatter}
	
\title{On the stability of fractional order Leslie-Gower type model with non-monotone functional response of intermediate predator}
\author[rmv]{Shuvojit Mondal\corref{cor1}}
\ead{shuvojitmondal91@gmail.com}
\author[cmbe]{Nandadulal Bairagi\fnref{fn2}}
\cortext[cor1]{Corresponding author}
\address[rmv]{Department of Mathematics, Rabindra Mahavidyalaya\\ Hooghly-712401, India.}
\address[cmbe]{Centre for Mathematical Biology and Ecology \\ Department of Mathematics, Jadavpur University\\ Kolkata-700032, India.}

\begin{abstract}
In this paper, an attempt is made to understand the dynamics of a fractional order three species Leslie-Gower predator prey food chain model with simplified Holling type IV functional response by considering fractional derivative in Caputo Sense. First, we prove different mathematical results like existence, uniqueness, non-negativity and boundedness of the solutions of fractional order dynamical system. The dissipativeness of the solution of the FDE system is discussed. Further, we investigate the Local stability criteria of all feasible equilibrium points. Global stability of the interior equilibrium point have also been discussed here. Using realistic parameter values,  numerically it has been observed that the fractional order system shows more complex dynamics, like chaos as fractional order becomes larger. Analytical results are illustrated with several examples in numerical section. 
\end{abstract}


\end{frontmatter}

\section{Introduction} \label{sec:1}

\setcounter{section}{1} \setcounter{equation}{0} 

Fractional calculus deals with the study of fractional-order integral and derivative operators over real or complex domains and their applications \cite{Rihan15}. In last two decades, fractional calculus have become more popular in various fields of science and engineering like physics, chemistry, biology and engineering \cite{Hashish12,Hilfer20,Laskin06,Rihan13}.  In ecological systems, fractional order models have also been used to understand the dynamics of interacting populations \cite{Ahmed07,RET13,CuiYang14,MET17,LET16,Vargas15,HET15,MET17a,MondalGaston19}. Lots of researchers are also showing interest to study the dynamics of discretized fractional-order systems and able to find more complex behaviors depending on both the step-size and fractional-order \cite{Mondal20,MondalCao20}. The major reason of using fractional derivatives is that it has the unique property of capturing the history of the variable, that is, it has memory \cite{Alidousti18,Matlob17}. In general, when the output of a system at each time $t$ depends only on the input at time $t$, such systems are usually known as memory-less systems. On the other side, when the system has to remember previous values of the input in order to determine the current value of the output, such systems are known as non-memory less systems, or memory systems \cite{Matlob17,Herrmann11}. Researchers consider the effect of recent memory as more important than the effect of older memory \cite{Caputo14} and the older memory can not be obtained by the help of integer order derivatives \cite{Baleanu12,Das11}. Due to the presence of memory, fractional derivative find wider area of applications and the models involving fractional derivatives provide a better agreement with real data than integer order derivative models. The value of the index of the fractional derivative, $m$ (say), characterizes the way in which the memory along different parts of the interval of integration affects the solution at a given time and it can be varied to best fit the real data \cite{Baleanu12,Diethelm20}. Therefore, it has been successfully observed the influence of memory concepts on the dynamics of systems \cite{Herrmann11,Butzer00}, and has been recently used in epidemiological models \cite{Saeedian17,Saka13}. 

After the seminal work of May \cite{May76}, exploring the chaotic behaviors in population models became a fascination. Lot of mathematical models have been proposed on food-chain and analyzed to show complicated dynamics like chaos \cite{Sahoo14,Huang01,Naji02,Naji10,Gakkhar12,Hasting91,Ali13}. Aziz-Alaoui \cite{Aziz02} discussed the complex dynamics in a modified Leslie-Gower three species food chain model with Holling type II response function. In \cite{Naji10}, Naji et al. also discussed about the chaotic dynamics in a modified Leslie-Gower food chain model with Bendigton-DeAnglis functional response. Gakkhar with Priyadarshi \cite{Gakkhar12} proposed a Leslie-Gower food cahin model. The functional response of Holling type IV describes a situation in which the predator’s per capita rate of predation decreases at sufficiently high prey densities \cite{Raw11}. Upadhyay in \cite{Raw11} used Holling type IV functional response of the form $\frac{wx}{d+x+x^2/i}$ to investigate the existence of complex dynamics in a three species food chain model. For details about this functional response, readers are referred to see \cite{Holling65}. Sokol and Howell in \cite{Sokol87} suggested a simplified Holling type IV function of the form $\frac{wx}{d+x^2}$ and found that it is simpler and better than the original function of Holling type IV. The more interesting fact is that they have all studied the ineteger order three species food chain models with this type of functional response.\\

In recent past, Alidousti and Ghahfarokhi \cite{Alidousti18} extended the work of Aziz-Alaoui \cite{Aziz02} and analyzed the following fractional order tri-trophic food chain model with Holling type II functional response:
\begin{eqnarray}\label{Tritophic fractional model}\nonumber
\frac{d^\alpha X}{dT^\alpha} & = &  X - \frac{\epsilon}{\delta} X^2 - \frac{v_0 XY}{\delta (d_0 + X)},~~ X(0)\geq 0, \\
\frac{d^\alpha Y}{dT^\alpha}& = & -\frac{g}{\delta} Y + \frac{v_1 XY}{\delta (d_1 + X)} - \frac{v_2 YZ}{\delta (d_2 + Y)},~~ Y(0)\geq 0, \\
\frac{d^\alpha Z}{dT^\alpha}& = & \frac{f}{\delta} Z^2 - \frac{v_3 Z^2}{\delta (d_3 + Y)},~~ Z(0)\geq 0, \nonumber
\end{eqnarray}
where $X, Y, Z$ are, respectively, the densities of prey, intermediate predator and top predator at any instant of at time $T$ and  $\alpha \in (0,1)$ is the order of the derivative. Sambath in \cite{Sambath18} also studied the asymptotic behavior of a fractional order three species predator-prey model with the same functional response. Investigations in fractional order Leslie-Gower type model with Holling type IV functional response is relatively less studied in population ecology. Therefore, in this paper, we consider a three species food chain model with simplified Holling type IV functional response to understand underlying dynamics of the model with respect to fractional order. \\

Recently, Ali et. al. \cite{Ali16} studied the following three-dimension coupled nonlinear autonomous system of integer order differential equations with non-monotone functional response (also called simplified Holling type IV functional response) to understand the underlying dynamics of food chain model:

\begin{eqnarray}\label{Tritophic_model}\nonumber
\frac{dX}{dT} & = & a_0 X - b_0 X^2 - \frac{v_0 XY}{d_1 + X^2},~~ X(0)\geq 0, \\
\frac{dY}{dT}& = & -a_1 Y + \frac{v_1 XY}{d_1 + X^2} - \frac{v_2 YZ}{d_2 + Y},~~ Y(0)\geq 0, \\
\frac{dZ}{dT}& = & Z\bigg(c_3 Z - \frac{v_3 Z}{d_3 + Y}\bigg),~~ Z(0)\geq 0, \nonumber
\end{eqnarray}
where $X, Y, Z$ are, respectively, the densities of prey, intermediate predator and top predator at any instant of at time $T$. This model considers interactions between a generalist top predator, specialist middle predator, and prey. Here, the specialist middle predator is consumed by the top predator, at a Holling type II rate. The interactions between the specialist middle predator and prey are modeled via a modified Holling type IV functional response. The interaction between the generalist top predator and specialist middle predator follow a modified Leslie–Gower scheme. That is the generalist top predator grows quadratically, because of sexual reproduction as $c_3Z^2$, and loses because of intra-species competition as $-\frac{v_3 Z^2}{d_3 + Y}$. The $d_3$ signifies that $Z$ is a generalist. The biological interpretation of all the parameters are described in the following table.  \\

{\small
	\begin{tabular}
		{|l l|} \hline
		& Table 1: Parameter interpretation \\ \hline
		Symbol & Meaning  \\  \hline
		
		$a_0  $                 & Growth rate of prey   \\ \hline
		$b_0 $                 & Intra specific competition coefficient \\ \hline
		$v_{i's} $                 & Maximum values that per-capita rate can attain \\	\hline
		$d_1  $                 & Measure of protection level provided by the environment to the prey \\ \hline
		$a_1 $                & Death rate of intermidiate predator \\	\hline
		$d_2  $                & Half-Saturation constant \\	\hline
		$c_3  $                & Growth rate of top predator via sexual reproduction \\	\hline
		$d_3  $                & Residual loss of top predator due to severe scarcity of it's favorite prey, Y \\	\hline
		
\end{tabular}}\\


All parameters are non-zero positive. \\


\noindent With the transformations
$$X = \frac{a_0}{b_0}x,~ Y = \frac{a_0^2}{b_0 v_0}y,~ Z = \frac{a_0^3}{b_0 v_0 v_2}z,~ T = \frac{t}{a_0},$$
the system (\ref{Tritophic_model}) takes the simplified form
\begin{eqnarray}\label{Tritophic_1_model}\nonumber
\frac{dx}{dt} & = & x(1-x)- \frac{sxy}{x^2 + a},~ x(0) = x_0\geq 0, \\
\frac{dy}{dt}& = & \frac{cxy}{x^2 + a} -by -  \frac{yz}{y + d},~ y(0) = y_0\geq 0, \\
\frac{dz}{dt}& = & pz^2 - \frac{qz^2}{y +r},~ z(0) = z_0\geq 0. \nonumber
\end{eqnarray}
where $$a = \frac{b_0^2 d_1}{a_0^2},~ b = \frac{a_1}{a_0},~ c = \frac{b_0 v_1}{a_0^2},~ d = \frac{d_2 v_0 b_0}{a_0^2},~ p = \frac{c_3 a_0^2}{b_0 v_0 v_2},~ q = \frac{v_3}{v_2},~ r = \frac{d_3 v_0 b_0}{a_0^2},~ s = \frac{b_0}{a_0}.$$

%

Starting from the integer-order three species Leslie~Gower type food chain model presented by (\ref{Tritophic_model}), we introduce the Caputo-type fractional order derivatives by replacing the usual integer-order derivatives to obtain the following fractional order system:
\begin{eqnarray}\label{Tritophic fractional order model}\nonumber
^{c}_{0} D^{m}_{T}X & = & a_0 X  - b_0 X^2 - \frac{v_0 XY}{d_0 + X^2}, \\
^{c}_{0} D^{m}_{T}Y & =  & -a_1 Y  + \frac{v_1 XY}{d_1 + X^2} - \frac{v_2 YZ}{d_2 + Y}, \\
^{c}_{0} D^{m}_{T}Z & =  & c_3 Z^2 - \frac{v_3 Z^2}{d_3 + Y}, \nonumber
\end{eqnarray}
with the initial conditions $X(0)\geq 0, Y(0) \geq 0, Z(0) \geq 0$, where $^{c}_{0} D^{m}_{t}$ is the Caputo fractional derivative \cite{Podlubny99} with fractional order $m$ $(0< m \leq1)$.  With the same transformations as before, the system (\ref{Tritophic fractional order model}) takes the following simplified form :

\begin{eqnarray}\label{Tritophic fractional order model_2}\nonumber
^{c}_{0} D^{m}_{t}x & = & x(1-x)- \frac{sxy}{x^2 + a},~ x(0) = x_0\geq 0, \\
^{c}_{0} D^{m}_{t}y & =  & \frac{cxy}{x^2 + a} -by -  \frac{yz}{y + d},~ y(0) = y_0\geq 0, \\
^{c}_{0} D^{m}_{t}z & =  & pz^2 - \frac{qz^2}{y +r},~ z(0) = z_0\geq 0, \nonumber
\end{eqnarray}
The state space of the system (\ref{Tritophic fractional order model_2}) is the positive cone $R^3_+ = \{(x, y, z) \in R^3: x\geq0, y\geq0, z\geq0 \}$. \\

\noindent We here intend to study the fractional order system (\ref{Tritophic fractional order model_2}), which is the simplified form of system (\ref{Tritophic fractional order model}) corresponding to integer order system (\ref{Tritophic_model}). We have shown that the solutions of system (\ref{Tritophic fractional order model_2}) are positively invariant and uniformly bounded in $R^3_+$ under some restrictions. Local stability criteria of the different equilibrium points have been discussed for fractional order system. Global stability of the interior equilibrium point have been only discussed here. Using realistic parameter values,  numerically it has been observed that the fractional order system shows more complex dynamics, like chaos as fractional order becomes larger and shows more simpler dynamics as the order $m$ decreases. Specially, it becomes stable for lower value of $m$. Simulation results are also given to validate the analytical results.

\section{Well-posedness}
Considering the biological significance of the model, we are only interested in solutions that are non-negative and bounded in the region $R^3_+ = \{(x, y, z) \in R^3: x\geq0, y\geq0, z\geq0 \}$. To prove the non-negativity and  uniform boundedness of our system, we shall use the following lemmas.\\

\noindent \textbf{Lemma 1}  \label{MVT} \cite{OdibatShawagfeh07} Suppose that$f(t)\in C[a,b]$ and $D^{m}_{a}f(t)\in C(a,b]$ with $0<m\leq1$. The Generalized Mean Value Theorem states that
	\begin{equation}\nonumber
	f(t) = f(a) + \frac{1}{\Gamma(m)}(D^{m}_{a}f)(\xi). (t-a)^m,
	\end{equation}
	where $a\leq\xi\leq t$, $\forall t \in (a,b]$.\\

From this lemma, one can easily prove the following result.\\

\noindent \textbf{Corollary 1} \cite{OdibatShawagfeh07} Suppose $f(t)\in C[a,b]$ and $^{c}_{t_0} D^{m}_{t}f(t)\in C(a,b)$, $0<m \leq 1$. If $^{c}_{t_0} D^{m}_{t}f(t)\geq 0, \forall t\in (a,b)$ then $f(t)$ is a non-decreasing function for each $t\in[a,b]$ and if $^{c}_{t_0} D^{m}_{t}f(t)\leq 0, \forall t\in (a,b)$ then $f(t)$ is a non-increasing function for each $t\in[a,b]$.\\

\noindent \textbf{Lemma 2}  \label{MVT_1}\cite{Hong-Li16} Let $u(t)$ be a continuous function on $[t_{0},\infty)$ and satisfying
	$$^{c}_{t_0} D^{m}_{t}u(t) \leq -\lambda u(t) + \mu,$$
	$$u(t_{0})                         =  u_{t_{0}}, $$
	where $0< m\leq1$, $(\lambda,\mu)\in\Re^{2}$, $\lambda\neq0$ and $t_{0}\geq0$ is the initial time. Then its solution has the form
	\begin{equation}\nonumber
	u(t) \leq \bigg(u_{t_{0}}- \frac{\mu}{\lambda}\bigg)E_{m}[-\lambda(t-t_{0})^{m}] + \frac{\mu}{\lambda}.\\
	\end{equation}

We have the following existence results on initial value problem (IVP) with caputo type fractional order differential equations.\\

\noindent \textbf{Lemma 3} \label{IVP_1}\cite{Elif12} Consider the initial value problem (IVP) with caputo type FDE
	\begin{equation}\label{IVP}
	^{c}_{t_{0}} D^{m}_{t} x(t) = f(t,x(t)),~~x(0) = x_0
	\end{equation}
	where $f \in C([0,T] \times R, R), ~~0<m<1$. Assume that $f \in C(R_0, R)$, where $R_0 = [(t,x) : 0\leq t\leq a, \mid x-x_0\mid \leq b]$ and let $\mid f(t,x)\mid \leq N$ on $R_0$. Then there exists at least one solution for the IVP (\ref{IVP}) on $0\leq t\leq \gamma$ where $\gamma = min \bigg(a, [~\frac{b}{M}\Gamma{(m +1)}]^{\frac{1}{m}}\bigg)$, $0<m < 1$.\\

\noindent \textbf{Lemma 4} \label{IVP_2} \cite{Elif12} Consider the initial value problem (IVP) given by (\ref{IVP}). Let \begin{equation}\nonumber
	g(v, x_*(v)) =  f(t -(t^m - v\Gamma{(m+1)})^{\frac{1}{m}}, ~x(t -(t^m - v\Gamma{(m+1)})^{\frac{1}{m}}))
	\end{equation}
	and assume that the conditions of lemma (\ref{IVP_1}) hold. Then, a solution of (\ref{IVP}), is given by
	\begin{equation}\nonumber
	x(t) = x_*(\frac{t^m}{\Gamma{(m+1)}})
	\end{equation}
	where $x_*(v)$ is the solution of the integer order differential equation
	\begin{equation}\nonumber
	\frac{d(x_*(v))}{dv} = g(v, x_*(v))
	\end{equation}
	with the initial condition $x_*(0) = x_0$.\\

\noindent \textbf{Theorem 1}  \label{Nonnegativity} All the solutions of system (\ref{Tritophic fractional order model_2}) which start in $R^3_+$ are non negative.\\

\noindent \textbf{Proof}: First we show that the solutions $x(t) \in \Re^{3}_{+} $ are non negative if it starts with positive initial values. If not, then there exists a $t_{1}> 0$ such that
\begin{eqnarray} \label{Function}\nonumber
&x(t)&> 0,~~0\leq t< t_{1},\\
&x(t)& = 0, ~~ t = t_1, \\
&x(t_{1}^+)& < 0.\nonumber
\end{eqnarray} 
Using (\ref{Function}) in the first equation of (\ref{Tritophic fractional order model_2}), we have
\begin{equation}
^{c}_{0} D^{m}_{t} x(t)|_{t = t_{1}} = 0.
\end{equation}
According to Lemma (\ref{MVT}), we have $x(t_{1}^{+}) = 0$, which contradicts the fact $x(t_{1}^{+}) < 0$. Therefore, we have $x(t)\geq0,~\forall ~t\geq0$. Using similar arguments, we can prove $y(t)\geq0, \forall t\geq0$ and $z(t)\geq0, \forall t\geq0$. So, it is proven that all the solutions of system (\ref{Tritophic fractional order model}) which start in $R^3_+$ are non negative. \\

\noindent Next we will show that, under some assumptions, all solutions $x(t), y(t)$ and $z(t)$ of our system (\ref{Tritophic fractional order model_2}) are uniformly bounded for sufficiently large $t$.\\

\noindent \textbf{Theorem 2} All the non negative solutions of system (\ref{Tritophic fractional order model_2}) which are initiating in $R^3_+$ are uniformly bounded, provided
	\begin{equation}\label{cond_1}
	\beta + \frac{\beta}{4b} +r < \frac{q}{p}
	\end{equation}
	and ultimately entering the region
	$$\Omega = \bigg\{(x,y,z) \in R^3_+ : 0\leq x \leq1, 0 \leq x+ \frac{y}{\beta} \leq 1 + \frac{1}{4b}, 0 \leq x + \frac{y}{\beta} + \alpha z \leq 1+ \frac{1}{4b} + \frac{M}{b} \bigg\},$$\\
	where $$\beta = \frac{v_1}{a_0}, ~\alpha = \frac{1}{b^2 (\beta +\frac{\beta}{4b} +r)}, ~M = \frac{1}{4(q-(\beta +\frac{\beta}{4b} +r)p)}.$$ Moreover, the system (\ref{Tritophic fractional order model_2}) is dissipative in $R^3_+$. \\

\noindent \textbf{Proof:} $(i)$ Let $(x(0), y(0), z(0)) \in \Omega$ and from the following theorem (\ref{Nonnegativity}), $(x(t)$ $, y(t), z(t))$ remain non negative in $R^3_+$, then we will show that $(x(t), y(t), z(t)) \in \Omega, \forall t \geq 0 $ and $\forall m \in (0,1]$. Now to reach our main objective, we have to prove the following steps for all $t \geq 0$ and $\forall m \in (0,1]$, \\
\textit{Step~(i-a):} $ x(t) \leq 1$;\\
\textit{Step~(i-b):} $ x(t) + \frac{y(t)}{\beta} \leq 1 + \frac{1}{4b}$;\\
\textit{Step~(i-c):} $ x(t) + \frac{y(t)}{\beta} + \alpha z(t) \leq 1+ \frac{1}{4b} + \frac{M}{b}$.\\

\textbf{\textit{Proof~of~Step~(i-a):}} We first prove that $ x(t) \leq 1, \forall t \geq 0$ and $\forall m \in (0,1]$. Since $x \geq 0, y\geq0, z\geq0$ in $R^3_+$, then any solution $\phi(t)  = (x(t), y(t), z(t))$ of (\ref{Tritophic fractional order model_2}), which starts in $R^3_+$, must satisfy the fractional order differential inequation
\begin{equation}\label{logistic}
^{c}_{0} D^{m}_{t}x  \leq  x(1-x), x(0) = x_0 > 0,  m \in (0, 1],\\
\end{equation}
which is clearly obtained from the first equation of (\ref{Tritophic fractional order model_2}). Moreover, this equation (\ref{logistic}) represents a fractional order logistic equation. Now we can apply lemma (\ref{IVP_2}) to solve this fractional order differential inequation  (\ref{logistic}). Here
\begin{equation}  \nonumber
g(v, x_*(v)) = x_*(v) (1-x_*(v))
\end{equation}
Then corresponding integer order differential inequation of this fractional IVP (\ref{logistic}) is
\begin{equation}  \nonumber
\frac{d(x_*(v))}{dv} \leq x_*(v) (1-x_*(v)), ~~x_*(0) = x_0.
\end{equation}
The solution of this integer order linear IVP is
\begin{equation} \nonumber
x_*(v) \leq \frac{1}{1+c_1 e^{-v}},\\
\end{equation}
where $c_1 = \frac{1}{x_0} - 1$. Consequently, the solution of the given fractional order IVP (\ref{logistic})  is
$$x(t) = x_*(\frac{t^m}{\Gamma{(m+1)}}) \leq \frac{1}{1+c_1 e^{-\frac{t^m}{\Gamma{(m+1)}}}} \leq 1, \forall t\geq 0, ~ 0<m \leq 1.$$

Then any solution $x(t)$ of (\ref{logistic}) must be bounded by $1$ with respect to any fractional order $m \in (0,1]$. Therefore, it follows that any non negative solution $\phi(t)$ of (\ref{Tritophic fractional order model_2}) satisfies $ x(t) \leq 1, \forall t \geq 0$ and $\forall m \in (0,1]$.\\

\textbf{\textit{Proof~of~Step~(i-b):}} We now prove that $ x(t) + \frac{y(t)}{\beta} \leq 1 + \frac{1}{4b}, \forall t \geq 0$ and $\forall m \in (0,1]$. Let us define a function
\begin{equation}
V_1(t) = x(t) + \frac{y(t)}{\beta},
\end{equation}
Taking fractional time derivative, we have
\begin{equation}\nonumber
^{c}_{0} D^{m}_{t} V_1(t) =~ ^{c}_{0} D^{m}_{t} x(t)+~ ^{c}_{0} D^{m}_{t} \frac{y(t)}{\beta} ~=~ x(1-x)- \frac{b}{\beta}y - \frac{1}{\beta}\frac{yz}{y + d}.
\end{equation}
Since all parameters are positive and solutions initiating in $R^3_+$ then,
\begin{equation}\nonumber
\begin{split}
^{c}_{0} D^{m}_{t} V_1(t) \leq & x(1-x)- \frac{b}{\beta}y,  \\
^{c}_{0} D^{m}_{t} V_1(t) \leq &  x(1-x) + bx - b (x+\frac{y}{\beta}), \\
^{c}_{0} D^{m}_{t} V_1(t) + & b V_1(t) \leq  b + \frac{1}{4},
\end{split}
\end{equation}
since in $\Omega$, $0\leq x\leq 1$ and $max_{[0, 1]} x(1-x)= \frac{1}{4}$. Applying Lemma (\ref{MVT_1}), we have
\begin{equation}\label{2D Solution}
\begin{split}
V_1(t) \leq & ~(V_1(0) - (1 + \frac{1}{4b}))E_{m}[-b t^{m}] + (1 + \frac{1}{4b}),\\
= & ~V_1(0) E_{m}[-b t^{m}] + (1 + \frac{1}{4b})(1 - E_{m}[-b t^{m}]).
\end{split}
\end{equation}
For $t\rightarrow \infty$, we thus have $V_1(t)\rightarrow (1 + \frac{1}{4b})$. Therefore, $V_1(t)\leq (1 + \frac{1}{4b}), \forall t \geq 0$ and $\forall m \in (0,1]$. Hence it follows that any non negative solution of (\ref{Tritophic fractional order model}) satisfies $ x(t) + \frac{y(t)}{\beta} \leq 1 + \frac{1}{4b}, \forall t \geq 0$ and $\forall m \in (0,1]$.\\

\textbf{\textit{Proof~of~Step~(i-c):}} We finally prove that $ x(t) + \frac{y(t)}{\beta} + \alpha z(t) \leq 1+ \frac{1}{4b} + \frac{M}{b}, \forall t \geq 0$ and $\forall m \in (0,1]$ holds, with $$\alpha = \frac{1}{b^2 (\beta +\frac{\beta}{4b} +r)}, ~M = \frac{1}{4(q-(\beta +\frac{\beta}{4b} +r)p)},$$ provided $\beta + \frac{\beta}{4b} +r < \frac{q}{p}$ . Again we define a function
\begin{equation}
V_2(t) = x(t) + \frac{y(t)}{\beta} + \alpha z(t),
\end{equation}
Taking fractional time derivative, we have
\begin{equation}\nonumber
\begin{split}
^{c}_{0} D^{m}_{t} V_2(t) = & ^{c}_{0} D^{m}_{t} x(t)+ ^{c}_{0} D^{m}_{t} \frac{y(t)}{\beta} + ^{c}_{0} D^{m}_{t} \alpha z(t),\\
= &  x(1-x)- \frac{b}{\beta}y - \frac{1}{\beta}\frac{yz}{y + d} + \alpha (p-\frac{q}{y+r}) z^2.
\end{split}
\end{equation}
Similarly to the previous step (i-b), since all parameters are positive, all solutions initiating in $Int(R^3_+)$ remain non negative and in $\Omega$, $0\leq x\leq 1$, $max_{[0, 1]} x(1-x)= \frac{1}{4}, y \leq \beta +\frac{\beta}{4b}$,  we get
\begin{eqnarray}\label{boundedness}\nonumber
^{c}_{0} D^{m}_{t} V_2(t) \leq & \frac{1}{4} + b - b V_2(t) + \alpha bz + \alpha (p-\frac{q}{y+r}) z^2,  \nonumber \\
^{c}_{0} D^{m}_{t} V_2(t) \leq & \frac{1}{4} + b - b V_2(t) + \alpha bz + \alpha (p-\frac{q}{\beta +\frac{\beta}{4b}+r}) z^2, \nonumber \\
^{c}_{0} D^{m}_{t} V_2(t) + & b V_2(t) \leq  b + \frac{1}{4} + M,
\end{eqnarray}

where, $$M = max_{z \in R^+} \bigg(\alpha bz + \alpha (p-\frac{q}{\beta +\frac{\beta}{4b}+r} )z^2\bigg).$$\\

Now we intend to find $M$, the maxima for the function $f(z) = \alpha bz + \alpha (p-\frac{q}{\beta +\frac{\beta}{4b}+r} )z^2, z\in R^+$. Here $f'(z) = \alpha b + 2 \alpha z (p-\frac{q}{\beta +\frac{\beta}{4b}+r} )$ and $f''(z) = 2 \alpha (p-\frac{q}{\beta +\frac{\beta}{4b}+r})$. Since $\beta + \frac{\beta}{4b} +r < \frac{q}{p}$, I observe $f''(z) < 0$ and hence $max [f(z)]$ exists at $$z = \frac{b}{2 (\frac{q}{\beta +\frac{\beta}{4b} +r} - p)} =z_1 (say).$$
Therefore using $\alpha = \frac{1}{b^2 (\beta +\frac{\beta}{4b} +r)}$, we have
\begin{equation}
M = max [f(z)]\mid_{ z = z_1} = \frac{\alpha b^2 (\beta +\frac{\beta}{4b} +r)}{4 (q - (\beta +\frac{\beta}{4b} +r)p)} = \frac{1}{4 (q - (\beta +\frac{\beta}{4b} +r)p)}.
\end{equation}

\noindent Then applying Lemma (\ref{MVT_1}) on (\ref{boundedness}), we have
\begin{equation}\label{3D Solution}
\begin{split}
V_2(t) \leq & ~(V_2(0) - (1 + \frac{1}{4b} + \frac{M}{b}))E_{m}[-b t^{m}] + (1 + \frac{1}{4b} + \frac{M}{b}), \\
= & ~V_2(0) E_{m}[-b t^{m}] + (1 + \frac{1}{4b} + \frac{M}{b})(1 - E_{m}[-b t^{m}]).
\end{split}
\end{equation}
For $t\rightarrow \infty$, we thus have $V_2(t)\rightarrow (1 + \frac{1}{4b} + \frac{M}{b})$. Therefore, $V_2(t)\leq (1 + \frac{1}{4b} + \frac{M}{b}), \forall t \geq 0$ and $\forall m \in (0,1]$. Hence it follows that any non negative solution of (\ref{Tritophic fractional order model}) satisfies for all $ m \in (0,1]$,
$$ x(t) + \frac{y(t)}{\beta} + \frac{1}{b^2 (\beta +\frac{\beta}{4b} +r)} z(t) \leq 1 + \frac{1}{4b} + \frac{1}{b} \frac{1}{4(q-(\beta +\frac{\beta}{4b} +r)p)}, \forall t \geq 0.$$ \\

Therefore, all the non negative solutions $x(t), y(t), z(t)$ of system (\ref{Tritophic fractional order model_2}) initiating in $R^3_+$ are uniformly bounded and entering the set $\Omega$, $\forall t \geq 0 $ and $\forall m \in (0,1]$. \\

From biological point of view, dissipativeness means all populations are bounded above. To prove this, we have to calculate the supremum of $x(t), V_1(t) = x(t) + \frac{y(t)}{\beta}$ and $V_2(t) = x(t) + \frac{y(t)}{\beta} + \alpha z(t)$ as $t \rightarrow +\infty$ for $0<m\leq 1$. The steps are following:

\noindent \textit{Step~(ii-a):} $ lim_{t \rightarrow +\infty} sup~x(t) \leq 1$;\\
\textit{Step~(ii-b):} $ lim_{t \rightarrow +\infty} sup~V_1(t) \leq 1 + \frac{1}{4b}$;\\
\textit{Step~(ii-c):} $ lim_{t \rightarrow +\infty} sup~V_2 \leq 1+ \frac{1}{4b} + \frac{M}{b}$.\\

\textbf{\textit{Proof~of~Step~(ii-a):}} Following theorem \ref{Nonnegativity}, since any non-negative solution $(x(t), y(t), z(t))$ of (\ref{Tritophic fractional order model_2}) satisfies $x(t)\leq 1, \forall t \geq 0$, so clearly for any $m \in (0, 1]$,  $lim_{t \rightarrow +\infty} sup~x(t) \leq 1$.\\

\textbf{\textit{Proof~of~Step~(ii-b):}} For the calculation of $ lim_{t \rightarrow +\infty} sup~V_1(t); V_1(t) =  x(t) + \frac{y(t)}{\beta}$, let $\epsilon > 0$ be given. Then there exists a $t_1 > 0$ such that $x(t) \leq 1 + \frac{\epsilon}{2}$ for all $t \geq t_1$. Applying $0 \leq x \leq 1$ and $max_{[0, 1]} x(1-x)= \frac{1}{4}$ for all $t \geq t_1 \geq 0$,  equation (\ref{2D Solution}) gives

\begin{equation}\label{2D Enter} \nonumber
\begin{split}
V_1(t) \leq & ~ 1 + \frac{1}{4b} - \bigg[1 + \frac{1}{4b} - \bigg(x(t_1) + \frac{y(t_1)}{\beta}\bigg)\bigg] E_{m}[-b (t - t_1)^{m}], \\
= & 1 + \frac{1}{4b} - \bigg[1 + \frac{1}{4b} - \bigg(x(t_1) + \frac{y(t_1)}{\beta}\bigg)\bigg] e^{-b (t - t_1)^{m}}, \\
\leq & 1 + \frac{1}{4b} - \bigg[1 + \frac{1}{4b} - \bigg(x(t_1) + \frac{y(t_1)}{\beta}\bigg)\bigg] e^{-b (t - t_1)}, [\because~0 < m \leq 1], \\
= & 1 + \frac{1}{4b} - \bigg[\bigg(1 + \frac{1}{4b}\bigg) e^{bt_1} - \bigg(x(t_1) + \frac{y(t_1)}{\beta}\bigg)e^{bt_1}\bigg] e^{-bt}, \\
\leq & 1 + \frac{1}{4b} - \bigg[\bigg(1 + \frac{1}{4b}\bigg) - \bigg(x(t_1) + \frac{y(t_1)}{\beta}\bigg)e^{bt_1}\bigg] e^{-bt},  \\
\leq & \bigg(1 + \frac{1}{4b} + \frac{\epsilon}{2}\bigg) - \bigg[\bigg(1 + \frac{1}{4b} +  \frac{\epsilon}{2}\bigg) - \bigg(x(t_1) + \frac{y(t_1)}{\beta}\bigg)e^{bt_1}\bigg] e^{-bt}, 
\end{split}
\end{equation}
for all  $t \geq t_1$. Suppose  $t_2 \geq t_1$ be such that $\bigg|\bigg(1 + \frac{1}{4b} +  \frac{\epsilon}{2}\bigg) - \bigg(x(t_1) + \frac{y(t_1)}{\beta}\bigg)e^{bt_1} \bigg| e^{-bt} \leq \frac{\epsilon}{2}$ for all $t \geq t_2$. Then we get $V_1(t) \leq 1 + \frac{1}{4b} + \epsilon$ for all $t \geq t_2$. Hence, $lim_{t \rightarrow +\infty} sup~V_1(t) \leq 1 + \frac{1}{4b}$, for any $m \in (0, 1]$.\\

\textbf{\textit{Proof~of~Step~(ii-c):}} 
Similarly we can consider $\epsilon > 0$. Then there exists a $t_3 > 0$ such that $V_1(t) \leq 1 + \frac{1}{4b} + \frac{\epsilon}{2}$ for all $t \geq t_3$. Next considering the equation (\ref{3D Solution}), we get for all $t \geq t_3 \geq 0$,

\begin{equation}\label{3D Enter}\nonumber
\begin{split}
V_2(t) \leq & ~ 1 + \frac{1}{4b}  + \frac{M}{b}- \bigg[1 + \frac{1}{4b} + \frac{M}{b} - \bigg(x(t_3) + \frac{y(t_3)}{\beta} \\
& + \alpha z(t_3)\bigg)\bigg] E_{m}[-b (t - t_3)^{m}], \\
= & 1 + \frac{1}{4b} + \frac{M}{b} - \bigg[1 + \frac{1}{4b} +\frac{M}{b} - \bigg(x(t_3) + \frac{y(t_3)}{\beta} + \alpha z(t_3)\bigg)\bigg] e^{-b (t - t_3)^{m}}, \\
\leq & 1 + \frac{1}{4b} +  \frac{M}{b} - \bigg[1 + \frac{1}{4b} +  \frac{M}{b}- \bigg(x(t_3) + \frac{y(t_3)}{\beta}\\
& + \alpha z(t_3)\bigg)\bigg] e^{-b (t - t_3)}, [\because~0 < m \leq 1], \\
= & 1 + \frac{1}{4b}  +  \frac{M}{b} - \bigg[\bigg(1 + \frac{1}{4b} +  \frac{M}{b}\bigg) e^{bt_3} - \bigg(x(t_3) + \frac{y(t_3)}{\beta}\\
& + \alpha z(t_3)\bigg)e^{bt_3}\bigg] e^{-bt}, \\
\leq & 1 + \frac{1}{4b} + \frac{M}{b} - \bigg[\bigg(1 + \frac{1}{4b} +  \frac{M}{b}\bigg) - \bigg(x(t_3) + \frac{y(t_3)}{\beta} + \alpha z(t_3)\bigg)e^{bt_3}\bigg] e^{-bt}, \\
\leq & \bigg(1 + \frac{1}{4b} + \frac{M}{b} + \frac{\epsilon}{2}\bigg) - \bigg[\bigg(1 + \frac{1}{4b} + \frac{M}{b}+  \frac{\epsilon}{2}\bigg) - \bigg(x(t_3) + \frac{y(t_3)}{\beta}\\
& + \alpha z(t_3)\bigg)e^{bt_3}\bigg] e^{-bt}.
\end{split}
\end{equation}
Again we suppose  $t_4 \geq t_3$ be such that $\bigg|\bigg(1 + \frac{1}{4b} + \frac{M}{b}+  \frac{\epsilon}{2}\bigg) - \bigg(x(t_3) + \frac{y(t_3)}{\beta} + \alpha z(t_3)\bigg)e^{bt_3}\bigg| e^{-bt} \leq \frac{\epsilon}{2}$ for all $t \geq t_4$. Then we get $V_2(t) \leq 1 + \frac{1}{4b}  + \frac{M}{b}+ \epsilon$ for all $t \geq t_4$. Hence $lim_{t \rightarrow +\infty} sup~V_2(t) \leq 1 + \frac{1}{4b} + \frac{M}{b}$, for any $m \in (0, 1]$. Therefore, we can say that our system (\ref{Tritophic fractional order model_2}) is dissipative in $R^3_+$ for all $m \in (0, 1]$.\\

Here we also study the existence and uniqueness of the solution of our system (\ref{Tritophic fractional order model_2}). We have the following Lemma due to Li et al \cite{LiChen10}.\\

\noindent \textbf{Lemma 5} \label{Uniqueness}\cite{LiChen10} Consider the system
	\begin{equation}\nonumber
	^{c}_{t_{0}} D^{m}_{t} x(t) = f(t,x), t>t_{0}
	\end{equation}
	with initial condition $x_{t_{0}}$, where $0<m\leq1$, $f:[t_{0},\infty)\times A\rightarrow\Re^{n}$, $A\in\Re^{n}$. If $f(t,x)$ satisfies the locally Lipschitz condition with respect to $x$ then there exists a unique solution of the above system on $[t_{0},\infty)\times A$ .\\

\subsection{Existence and uniqueness}

Using Lemma (\ref{Uniqueness}), here we study the existence and uniqueness of the solution of system (\ref{Tritophic fractional order model_2}) in the region $A \times[0, T_1]$, where $A = \{(x,y,z)\in\Re^{3}|~ max\{|x|, |y|, |z|\} \leq M_1\}$, $T_1<\infty$ and $M_1$ is large. Denote $X = (x,y,z)$, $\bar{X} = (\bar{x}, \bar{y}, \bar{z})$. Consider a mapping $H: A \rightarrow\Re^{3}$ such that $H(X) = (H_{1}(X), H_{2}(X), H_{3}(X))  $, where
\begin{eqnarray}\label{Existence}
H_{1}(X) & = &  x(1-x)- \frac{sxy}{x^2 + a},\nonumber \\
H_{2}(X) & = & \frac{cxy}{x^2 + a} -by -  \frac{yz}{y + d}, \\ 
H_3(x) & =  & pz^2 - \frac{qz^2}{y +r}. \nonumber
\end{eqnarray}
For any $X, \bar{X} \in A$, it follows from (\ref{Existence}) that \\
\begin{eqnarray}\nonumber
\begin{split}
\parallel H(X) - H(\bar{X})\parallel = & \mid H_{1}(X)-H_{1}(\bar{X})\mid +\mid H_{2}(X)-H_{2}(\bar{X})\mid \\
& + \mid H_{3}(X)-H_{3}(\bar{X})\mid \\
= & \mid  x(1-x)- \frac{sxy}{x^2 + a} -  \bar{x}(1-\bar{x}) + \frac{s\bar{x}\bar{y}}{\bar{x}^2 + a}\mid \\
& + \mid \frac{cxy}{x^2 + a} -by -  \frac{yz}{y + d} - \frac{c\bar{x}\bar{y}}{\bar{x}^2 + a} + b\bar{y} + \frac{\bar{y}\bar{z}}{\bar{y} + d}\mid \\
& + \mid  pz^2 - \frac{qz^2}{y +r} -  p\bar{z}^2  + \frac{q\bar{z}^2}{\bar{y} +r}\mid \\
= & \mid (x-\bar{x}) - (x^{2} -\bar{x}^{2}) - s\bigg(\frac{xy}{x+a} - \frac{\bar{x}\bar{y}}{\bar{x} +a}\bigg)\mid \\
& +\mid c \bigg(\frac{xy}{x +a} - \frac{\bar{x}\bar{y}}{\bar{x} +a}\bigg) - b(y - \bar{y}) - (\frac{yz}{y+d} - \frac{\bar{y}\bar{z}}{\bar{y} + d})\mid \\
& + \mid p(z^{2} - \bar{z}^{2}) - q(\frac{z^2}{y+r} - \frac{\bar{z}^2}{\bar{y} +r})\mid \\
& \leq \mid x-\bar{x}\mid + \mid x^{2} -\bar{x}^{2}\mid \\
&+ (s+c)\mid\bigg(\frac{xy}{x +a} -\frac{\bar{x}\bar{y}}{\bar{x} +a}\bigg)\mid + b\mid y - \bar{y}\mid \\
\end{split}
\end{eqnarray}
\begin{eqnarray}\nonumber
\begin{split}
& +\mid \frac{yz}{y+d} - \frac{\bar{y}\bar{z}}{\bar{y} + d}\mid + p \mid z^{2} - \bar{z}^{2}\mid + q\mid (\frac{z^2}{y+r} - \frac{\bar{z}^2}{\bar{y} +r})\mid \\
& \leq \mid x-\bar{x}\mid + 2M_1 \mid x - \bar{x}\mid + \frac{(s+c)}{a^2}\mid a(xy - \bar{x}\bar{y})\\
& + x\bar{x}(\bar{x}y - \bar{y}x)\mid +b\mid y - \bar{y}\mid + \frac{1}{d^2}\mid d(yz - \bar{y}\bar{z}) + y \bar{y}(z-\bar{z})\mid \\
& + p \mid z^{2} - \bar{z}^{2}\mid + \frac{q}{r^2} \mid r(z-z^2) + (z^2\bar{y} - \bar{z}^2 y)\mid \\
& \leq \bigg(1 + 2M_1 + \frac{M_1(s+c)}{a} + \frac{M_1^3(s+c)}{a^2}\bigg)\mid x -\bar{x}\mid \\
&+ \bigg(M_1(\frac{s+c}{a} +\frac{1}{d}) +b + M_1^2 (\frac{M_1(s+c)}{a^2} + \frac{q}{r^2})\bigg)\mid y-\bar{y}\mid \\
& + \bigg(M_1(2p + \frac{1}{d} + \frac{2q}{r}) + M_1^2 (\frac{1}{d^2} + \frac{2q}{r^2})\bigg)\mid z-\bar{z}\mid\\ 
& \leq L\parallel (x,y,z) - (\bar{x}, \bar{y}, \bar{z})   \parallel \\
& \leq L\parallel X - \bar{X} \parallel,
\end{split}
\end{eqnarray}
where $L = max \{1 + 2M_1 + \frac{M_1(s+c)}{a} + \frac{M_1^3(s+c)}{a^2}, M_1(\frac{s+c}{a} +\frac{1}{d}) \\
+b + M_1^2 (\frac{M_1(s+c)}{a^2} + \frac{q}{r^2}), M_1(2p + \frac{1}{d} + \frac{2q}{r}) + M_1^2 (\frac{1}{d^2} + \frac{2q}{r^2})\}.$\\

Thus $H(X)$ satisfies Lipschitz condition with respect to $X$ and following Lemma (\ref{MVT_1}), there exists a unique solution $X(t)$ of the system (\ref{Tritophic fractional order model_2}) with the initial condition $X(0) = (x(0), y(0), z(0))$.

\section{Existence and stability of equilibria}
\noindent We have the following stability result on fractional order dynamical systems.\\

\noindent \textbf{Theorem 3} \cite{MET17a,Petras11} Consider the following fractional order system
	\begin{equation}\label{Stability condition} \nonumber
	^{c}_{0} D^{m}_{t} x(t) = f(x), x(0)= x_{0}
	\end{equation}
	with $0<m\leq1, x\in \Re^{n}$ and $f: \Re^{n}\rightarrow \Re^{n}$. The equilibrium points of the  above system are calculated by solving the equation: $f(x) = 0$. These equilibrium points are locally asymptotically stable if all eigenvalues $\lambda_{i}$ of the Jacobin matrix $J = \frac{\partial f}{\partial x}$ evaluated at the equilibrium points satisfy $\mid arg(\lambda_{i})\mid >\frac{m \pi}{2}$, $i = 1,2,----,n$.\\

\noindent An equilibrium point of system (\ref{Tritophic fractional order model_2}) is found by solving the three equations $D^{m}_{t} x(t) =  D^{m}_{t} y(t) = D^{m}_{t} z(t) = 0$ in (\ref{Tritophic fractional order model_2}). There are  four biologically feasible non-negative equilibrium points of system (\ref{Tritophic fractional order model_2}). The trivial equilibrium  $E_{0} = (0, 0, 0)$ and the axial equilibrium $E_{1} = (1,0,0)$ are always exist.  The planner equilibrium point $E_2 = (\bar{x}, \bar{y}, 0)$ exists uniquely in the positive quadrant of $xy-$ plane, where $\bar{x}$ and $\bar{y}$ are given by $$\bar{x} = \frac{c}{2b}, ~\bar{y} = \frac{1}{s}(1 - \bar{x})(a + \bar{x}^2),$$
provided that the following conditions are hold:
\begin{equation}\label{cond}
\frac{c}{2b} < 1, ~c^2 - 4ab^2 = 0. 
\end{equation}
We observe that in the absence of prey $x$, both predators $y$ and $z$ can not survive. So there is no equilibrium point in the $yz-$ plane. Similarly we can also conclude that there is no equlibrium point in $xz-$ plane. Now there exists a unique interior equilibrium point
$E^* = (x^*, y^*, z^*)$ of the system (3), where the equilibrium population densities are given by
\begin{eqnarray}\label{Equilibrium relation_1}
y^{*} = \frac{q}{p} - r,
\end{eqnarray}
while $x^*$ is the positive root of the cubic equation

\begin{eqnarray}\label{Equilibrium relation_2}
x^3 - x^2 + ax + (sy^* - a) = 0,
\end{eqnarray}
this equation can be written as 
\begin{eqnarray}\label{Equilibrium relation_3}
f(x) = Ax^3 + Bx^2 + Cx + D = 0,
\end{eqnarray}
where $A = 1, B = -1, C = a$ and $D = (sy^* - a)$. Now since $0 \leq x^* \leq 1$, then $f(0) = D <0,$ if $y^* < \frac{a}{s}$ and $f(1) = sy^* > 0$. Thus, $f(0)f(1) = sy^*(sy^* - a) < 0$, and then there is a positive root of equation (\ref{Equilibrium relation_3}) lies in $(0, 1)$ when $y^* < \frac{a}{s}$ is satisfied. Now from the second eqaution of system (\ref{Tritophic fractional order model_2}), we obtain
\begin{eqnarray}\label{Equilibrium relation_4}
~z^* = (-b + \frac{cx^*}{a +{x^*}^2})(y^* + d),
\end{eqnarray}
and it exists if $b < \frac{cx^*}{a +{x^*}^2}$. Therefore the positivity condition of $E^*$ in $R^3_+$ are $$y^* < \frac{a}{s}, ~b < \frac{cx^*}{a +{x^*}^2},$$ where $v_3 > c_3 d_3$. Different stability results for the equilibrium points $E_0, E_1$, $E_2$ and $E^*$ are given in the following.\\

\noindent Now to investigate the dynamical behavior of the equilibrium points $E_i, ~(i = 0,1,2)$ and $E^*$, we first construct the Jacobian matrix $J$ evaluated at an equilibrium point $(x, y, z)$ of the system (\ref{Tritophic fractional order model_2}) is
\begin{equation}\label{jacobian}
J(x, y, z) = \begin{pmatrix}
a_{11} & a_{12} & a_{13}\\a_{21} & a_{22} & a_{23} \\ a_{31} & a_{32} & a_{33}
\end{pmatrix},
\end{equation}
where $a_{11} = 1- 2x - \frac{sy(a - x^2)}{(a+x^2)^2},~~a_{12} = -\frac{sx}{a+x^2},~~ a_{13} = 0, ~~a_{21} = \frac{cy(a - x^2)}{(a+x^2)^2},\\
~~a_{22} = -b + \frac{cx}{a+x^2} - \frac{dz}{(d+y)^2}, a_{23} = -\frac{y}{d+y},~~a_{31} = 0, ~~ a_{32} = \frac{qz^2}{(r +y)^2},\\
~~a_{33} = 2z(p - \frac{q}{y+r}).$ \\

Then the Jacobian matrices evaluated at $E_0$, $E_1$ and $E_2$ are given by
\[
J(E_0) = \begin{pmatrix} 
1 & 0 & 0 \\ 
0 & -b   & 0 \\ 
0 & 0 & 0
\end{pmatrix},\]
\[
J(E_1) = \begin{pmatrix} 
-1 & -\frac{s}{a+1} & 0 \\ 
0 & -b + \frac{c}{a+1} & 0 \\ 
0 & 0 & 0
\end{pmatrix},\]
\[
J(E_2) = \begin{pmatrix} 
1 - 2\bar{x} - \frac{(1-\bar{x})(a - \bar{x}^2)}{a+\bar{x}^2}& -\frac{s\bar{x}}{a+\bar{x}^2} & 0 \\ 
\frac{c(1-\bar{x})(a - \bar{x}^2)}{s(a + \bar{x}^2)}& 0 & -\frac{\bar{y}}{(d + \bar{y})} \\
0 & 0 & 0
\end{pmatrix}.\]

Clearly, the eigenvalues of $J(E_0)$ are $\xi_{1} = 1$, $\xi_{2} = -b$ and $\xi_{3} = 0$. Note that $arg(\xi_{3})$ is undefined. Since one them is a positive real and anothe one is a negative real, then $E_0$ is always unstable. Therefore $E_{0}$ is non-hyperbolic. \\

Next, the eigenvalues of $J(E_1)$ are $\xi_{1} = -1 ~(<0)$, ~$\xi_{2} = \frac{c - b - ab}{1+a}$ and $\xi_{3} = 0$. Hence $E_{1}$ is also non-hyperbolic. Note that If $c - b > ab$ then $\xi_{2}>0$. In this case, $E_1$ is always unstable saddle along $x-$ direction. If $c - b < ab$ then $\xi_{2}<0$. Consequently, two of the eigenvalues are negative real, so in this case $E_1$ is stable manifold along $x$ and $y-$ direction.\\

Again from the variational matrix of $E_2$, the eigenvalues of $J(E_2)$  are $\xi_{1,2} = \frac{1}{2}[P \pm \sqrt{P^{2} - 4Q}]$, where $P = 1 - 2\bar{x} - \frac{(1-\bar{x})(a - \bar{x}^2)}{a+\bar{x}^2}, ~Q = \frac{c\bar{x}(1-\bar{x})(a- \bar{x}^2)}{(a + \bar{x}^2)^2}$ and  $\xi_{3} = 0$. Since one of the eigenvalue $\xi_3$ becomes zero, so $E_2$ is non-hyperbolic equilibrium point. \\


\noindent For local stability of the interior equilibrium $E^*$, we compute the Jacobian matrix of system (\ref{Tritophic fractional order model_2}) at $E^* = (x^*, y^*, z^*)$ as
\begin{equation}
J(E^{*}) = \begin{pmatrix}
1 - 2x^* - \frac{(1-x^*)(a - {x^*}^2)}{a+{x^*}^2}& -\frac{sx^*}{a+{x^*}^2} & 0 \\ \frac{c(1-x^*)(a - {x^*}^2)}{s(a + {x^*}^2)} & \frac{y^*z^*}{(y^*+d)^2} & -\frac{y^*}{d + y^*}  \\ 0 & \frac{p{z^*}^2}{y^* + r} & 0
\end{pmatrix}.
\end{equation}
The eigenvalues are the roots of the cubic equation
\begin{equation}\label{cubic equation}
F(\xi) = \xi^{3} + A_{1} \xi^{2} + A_{2} \xi + A_{3} =   0,
\end{equation}
where 
$
A_{1} = -1 + 2x^* + \frac{(1-x^*)(a - {x^*}^2)}{a+{x^*}^2} - \frac{y^*z^*}{(y^*+d)^2},\\~ 
A_{2} = \frac{y^*z^*}{(y^*+d)^2}\bigg(1 - 2x^* - \frac{(1-x^*)(a - {x^*}^2)}{a+{x^*}^2}\bigg) + \frac{cx^*(1-x^*)(a - {x^*}^2)}{(a + {x^*}^2)^2} + \frac{py^*{z^*}^2}{(y^* + d)(y^* + r)},\\
A_{3} = -\frac{py^*{z^*}^2}{(y^* + d)(y^* + r)}\bigg(1 - 2x^* - \frac{(1-x^*)(a - {x^*}^2)}{a+{x^*}^2}\bigg).$\\

\noindent The equilibrium $E^*$ is said to be locally asymptotically stable if all eigenvalues of (\ref{cubic equation}) satisfy  $\mid arg(\xi_{i})\mid >\frac{m \pi}{2}, \forall m \in (0,1]$, $i = 1,2,3$. One can then determine the stability of $E^*$ by noting the signs of the coefficients $A_i$ and discriminant $D(F)$ of the cubic polynomial $F(\xi)$ \cite{Ahmed07,Ahmed06}. The discriminant $D(F)$ of the cubic polynomial $F(\xi)$ is
\[
\mathbf{D(F)} = - \begin{vmatrix}
1 & A_{1} & A_{2} & A_{3} & 0 \\ 0 & 1 & A_{1} & A_{2} & A_{3} \\ 3 & 2A_{1} & A_{2} & 0 & 0 \\ 0 & 3 & 2A_{1} & A_{2} & 0 \\ 0 & 0 & 3 & 2A_{1} & A_{2}
\end{vmatrix}= 18A_{1}A_{2}A_{3} + (A_{1}A_{2})^2 - 4A_{3}A_{1}^{3} - 4A_{2}^{3} - 27A_{3}^{2}.
\]
Then the following theorem regarding local asymptotic stability of $E^*$ of the system (\ref{Tritophic fractional order model_2}) is true \cite{Ahmed07,MET17a,Ahmed06}.\\

\noindent \textbf{Theorem 4} \label{Interior Stability}
	\begin{itemize}
		\item[(i)] If $D(F) > 0$,  $ A_{1}>0$, $A_{3}>0$ and $A_{1}A_{2}- A_{3}>0$ then the interior equilibrium $E^{*}$ is locally asymptotically stable for all $m \in (0, 1]$.
		
		\item[(ii)] If $D(F) < 0$, $A_{1}\geq 0$, $A_{2}\geq 0$, $A_{3} > 0$ and $0< m < \frac{2}{3}$ then the interior equilibrium $E^{*}$ is locally asymptotically stable.
		
		\item[(iii)] If $D(F) < 0$, $A_{1} < 0$, $A_{2} < 0$ and $m > \frac{2}{3}$ then the interior equilibrium $E^{*}$ is unstable.
		\item[(iv)] If $D(F) < 0$, $A_{1} > 0$, $A_{2} > 0$, $A_{1}A_{2} = A_{3}$ and $0< m < 1$ then the interior equilibrium $E^{*}$ is locally asymptotically stable.
	\end{itemize}

\noindent \textbf{Proof:} ~(i) If $D(F)$ is positive then all the roots of (\ref{cubic equation}) are real and distinct. If not, let us assume that $F(\xi) = 0$ has one real root $\xi_1$ and another two complex conjugate roots $\xi_2$, $\xi_3$. In terms of the roots, the discriminant of $F(\xi)$ can be written as \cite{Janson07}
\begin{equation}\label{Discriminant eqn}
D(F) =[( \xi_1 - \xi_2)(\xi_1 - \xi_3)(\xi_2 - \xi_3)]^2.
\end{equation}
Note that
\begin{eqnarray}
\begin{split}
( \xi_1 - \xi_2)(\xi_1 - \xi_3)(\xi_2 - \xi_3) = & ( \xi_1 - \xi_2)(\xi_1 - \overline{\xi_2})(\xi_2 - \overline{\xi_2})\\
= & ( \xi_1 - \xi_2)(\xi_1 - \overline{\xi_2})2 Im(\xi_2)i\\
= & ( \xi_1 - \xi_2)\overline{(\xi_1 - \xi_2)}2 Im(\xi_2)i\\
= & 2|\xi_1 - \xi_2|^2 Im(\xi_2)i.
\end{split}
\end{eqnarray}
Thus, \begin{equation}
D(F) =[2|\xi_1 - \xi_2|^2 Im(\xi_2)i]^2 < 0,
\end{equation}
which contradicts the fact that $D(F) > 0$. Therefore, whenever $D(F)>0$ then $F(\xi) = 0$ has three real distinct roots. Since $ A_{1}>0$, $A_{3}>0$ and $A_{1}A_{2}- A_{3}>0$, all roots of $F(\xi) = 0$ has negative real roots or complex conjugate roots with negative real parts. As $D(F) > 0$, so all roots of $F(\xi) = 0$ are real negative. Consequently, $\mid arg(\xi_{i})\mid  = \pi > \frac{\alpha \pi}{2}, \forall \alpha\in (0,1]$, $i = 1, 2, 3$, and the equilibrium $E^{*}$ is  locally asymptotically stable. This completes the proof of (i).\\

\noindent (ii) We have seen in (i) that $F(\xi) = 0$ has one real and two complex conjugate roots if $D(F) < 0$. Since $A_{3} > 0$, following (\ref{cubic equation}), the real root is negative. We thus consider the roots as $\xi_{1} = -b$, $(b\in R_+)$ and $\xi_{2,3} = \beta \pm i \gamma,$ $(\beta, \gamma \in R)$ and
$$F(\xi) = (\xi +b)(\xi - \beta - i \gamma)(\xi - \beta + i \gamma).$$ Comparing this with (\ref{cubic equation}), we have
$A_{1} =  b-2\beta,
A_{2} = \beta^2 + \gamma^2 - 2b\beta,
A_{3} = b(\beta^2 + \gamma^2).$
Now $A_{1} \geq 0 \Rightarrow b \geq 2\beta$. Noting $\beta^2 \sec^2\theta= \beta^2 + \gamma^2$ and $A_{2} \geq 0$, we have $sec^2\theta \geq  4.$ Therefore, $\theta= |arg(\xi)| \geq \frac{\pi}{3}$. Since $0< \alpha < \frac{2}{3}$, then $ |arg(\xi)| = \theta \geq \frac{\pi}{3} > \frac{\alpha \pi}{2}$ holds. Thus, all roots of (\ref{cubic equation}) satisfy $\mid arg(\xi_{i})\mid >\frac{\alpha \pi}{2}, \forall \alpha\in (0,1]$ and the equilibrium $E^{*}$ is locally asymptotically stable. This completes the proof of (ii). Proof of $(iii)$ is similar to the proof of (ii) and hence omitted.\\

\noindent Since $D(F) < 0$, $A_{1} > 0$, $A_{2} > 0$, from the previous case, we have the
\begin{eqnarray}\nonumber
A_{1} =  b-2\beta, ~A_{2} =  \beta^2 + \gamma^2 - 2b\beta, ~A_{3} =  b(\beta^2 + \gamma^2).
\end{eqnarray}
Note that $A_{1} > 0 \Rightarrow b> 2\beta$, $A_{2} > 0 \Rightarrow \beta^2 + \gamma^2 - 2b\beta > 0$
and $A_{1}A_{2} = A_{3} \Rightarrow (b-2\beta)(\beta^2 + \gamma^2 - 2b\beta) = b(\beta^2 + \gamma^2) \Rightarrow \beta(b^2 + \beta^2 + \gamma^2 - 2b\beta) = 0$. Then two cases arise: \\
\noindent\textbf{Case 1:} If $\beta = 0$ then three roots $\xi_{1},\xi_{2},\xi_{3}$ of (\ref{cubic equation}) are $-b, \pm i\gamma$. One can see that $\mid arg(\xi_{1})\mid  = \pi > \frac{\alpha \pi}{2}$ and $\mid arg(\xi_{2,3})\mid  = \frac{\pi}{2} > \frac{\alpha \pi}{2}, \forall\alpha\in (0,1)$. Therefore, the equilibrium $E^{*}$ is locally asymptotically stable.\\
\noindent\textbf{Case 2:} If $b^2 + \beta^2 + \gamma^2 -2b\beta=0$ then we have $b=\beta$ and $\gamma=0$. Using it in $b> 2\beta$ and $\beta^2 + \gamma^2 >2b\beta$, we obtain $b<0$, which contradicts the assumption $b\in R_+$.\\
Thus, if $A_{1} > 0$, $A_{2} > 0$, $A_{1}A_{2} = A_{3}$ then one root is real negative and the other two are purely imaginary and therefore $\mid arg(\xi_{1})\mid  = \pi > \frac{\alpha \pi}{2}$ and $\mid arg(\xi_{2,3})\mid  = \pi/2 > \frac{\alpha \pi}{2},$ $\forall \alpha\in (0,1)$, implying local asymptotic stability of $E^{*}$. This completes the proof.\\

\noindent Then, we proceed to prove the global stability results of the interior equilibrium point  $E^* = (x^*, y^*, z^*)$. The following lemma will be used in proving it.\\

\noindent \textbf{Lemma 6}\label{Global} \cite{Vargas15} Let $x(t)\in \Re_{+}$ be a continuous and derivable function. Then for any time instant $t>t_{0}$
	\begin{equation}\nonumber
	^{c}_{t_{0}} D^{m}_{t}\bigg[x(t) - x^{*} - x^{*}ln\frac{x(t)}{x^{*}}\bigg] \leq \bigg(1-\frac{x^{*}}{x(t)}\bigg)~~{^{c}_{t_{0}}} D^{m}_{t}x(t), x^{*}\in \Re_{+}, \forall m\in(0,1].
	\end{equation}

\noindent \textbf{Theorem 5} \label{Global_Stability} The interior equilibrium $E^* = (x^*, y^*, z^*)$ of system (\ref{Tritophic fractional order model_2}) is globally asymptotically stable for any $m\in (0, 1]$ if
	\begin{eqnarray}\label{cond_2} \nonumber
	(i)&\frac{2sy^*}{a({x^*}^2 +a)} + \frac{s}{2a^2}- 1 <0, \nonumber \\
	(ii)&\frac{cx^* - b({x^*}^2 + a)}{a\beta d} + \frac{s}{2a^2}+ \frac{1}{2}\bigg(\frac{q}{4br\alpha(q - p(\beta+\frac{\beta}{4b}+r))} - \frac{{x^*}^2 +a}{a\beta(\beta+\frac{\beta}{4b}+d)}\bigg)  < 0, \nonumber \\
	(iii)&\frac{q}{4br\alpha(q - p(\beta+\frac{\beta}{4b}+r))} - \frac{{x^*}^2 +a}{a\beta(\beta+\frac{\beta}{4b}+d)} < 0,\nonumber
	\end{eqnarray}
	where $\alpha = \frac{1}{b^2 (\beta +\frac{\beta}{4b} +r)}>0.$

\noindent \textbf{Proof:} Let us consider the Lyapunov function
\begin{eqnarray}\nonumber
\begin{split}
V(x,y,z) &= \bigg(x - x^* - x^*ln\frac{x}{x^*}\bigg) + \frac{({x^*}^2 + a)}{a\beta}\bigg(y - y^* - y^*ln\frac{y}{y^*}\bigg) \\
&+ (y^* + r)\bigg(z - z^* - z^*ln\frac{z}{z^*}\bigg).
\end{split}
\end{eqnarray}
It is easy to see that $V = 0$ only at $(x, y, z) = (x^*, y^*, z^*)$ and $V > 0$ whenever $(x, y, z) \neq (x^*, y^*, z^*)$. Considering the $m-th$ order fractional derivative of $V(x,y,z)$ along the solutions of (\ref{Tritophic fractional order model_2}), we have
\begin{eqnarray}
\begin{split}
^c_{0}D^{m}_{t}V(x,y,z) &= {^c_{0}}D^{m}_{t}\bigg(x - x^* - x^*ln\frac{x}{x^*}\bigg) + \frac{({x^*}^2 + a)}{a\beta} {^c_{0}}D^{m}_{t}\bigg(y - y^* - y^*ln\frac{y}{y^*}\bigg) \\
&+ (y^* + r){^c_{0}}D^{m}_{t}\bigg(z - z^* - z^*ln\frac{z}{z^*}\bigg). \nonumber
\end{split}
\end{eqnarray}
Using Lemma (\ref{Global}) and making some algebraic manipulations, we have
\begin{eqnarray}\label{Globality_2}\nonumber
\begin{split}
^c_{0}D^{m}_{t}V(x,y,z) \leq & \frac{(x - x^*)}{x} {^c_{0}}D^{m}_{t}x(t) + \frac{({x^*}^2 + a)}{a\beta} \frac{(y - y^*)}{y} {^c_{0}}D^{m}_{t}y(t) \\
& + (y^* + r)\frac{(z - z^*)}{z} {^c_{0}}D^{m}_{t}z(t) \\
=& (x - x^*) \bigg[(1- x) - \frac{sy}{x^2 +a}\bigg] + \frac{({x^*}^2 + a)}{a}(y - y^*)\\
& \bigg[\frac{sx}{x^2 +a} - \frac{b}{\beta} - \frac{z}{\beta(d+y)}\bigg] + (y^* + r)(z-z^*)\bigg[pz - \frac{qz}{y+r}\bigg]\\
\end{split}
\end{eqnarray}
\begin{eqnarray}
\begin{split}
=& (x - x^*)\bigg[(x^* - x) + \frac{sy^*}{{x^*}^2 +a} - \frac{sy}{x^2 +a}\bigg] \\
&+ \frac{({x^*}^2 + a)}{a}(y - y^*)\bigg[\frac{sx}{x^2 + a} - \frac{sx^*}{{x^*}^2 + a} +\frac{z^*}{\beta(d+y^*)} - \frac{z}{\beta(d+y)}\bigg] \\
&+ z(y^* + r)(z-z^*)\bigg[\frac{q}{y^* +r} - \frac{q}{y+r}\bigg]\\
=&  - (x - x^*)^{2}  + s(x - x^*) \bigg[\frac{y^*(x^2 + a) -y({x^*}^2 + a)}{(x^2 +a)({x^*}^2 +a)}\bigg] \\
& + \frac{({x^*}^2 + a)}{a} s(y - y^*) \bigg[\frac{x({x^*}^2 + a) - x^*(x^2 + a)}{(x^2 +a)({x^*}^2 +a)}\bigg]\\
&+ \frac{({x^*}^2 + a)}{a} (y - y^*) \bigg[\frac{d(z^* -z) + \{y^*(z^* -z) + z^* (y - y^*)\}}{(y +d)(y^* +d)}\bigg] \\
& + \frac{qz}{y+r}(y-y^*)(z-z^*)\\
= & \bigg[\frac{sy^*(x+x^*)}{(x^2 +a)({x^*}^2 +a)} - 1\bigg] (x - x^*)^{2} - \frac{sxx^*}{a({x^*}^2 +a)}(x-x^*)(y-y^*) \\
&+ \frac{({x^*}^2 + a)}{a\beta}\frac{z^*(y - y^*)^2}{(y+d)(y^* + d)} \\
&+ \bigg[\frac{qz}{y+r} - \frac{{x^*}^2 + a}{a\beta(y+d)}\bigg](y-y^*)(z-z^*)\\
\leq & \bigg[\frac{2sy^*}{a({x^*}^2 +a)} - 1\bigg] (x - x^*)^{2}  + \frac{z^*({x^*}^2 + a)(y - y^*)^2}{a\beta d(y^* + d)} \\
&+\frac{sxx^*}{a({x^*}^2 +a)}\bigg[\frac{(x-x^*)^2 + (y-y^*)^2}{2}\bigg] \\
+& \bigg[\frac{q}{4br\alpha(q - p(\beta+\frac{\beta}{4b}+r))} - \frac{{x^*}^2 +a}{a\beta(\beta+\frac{\beta}{4b}+d)}\bigg]\bigg[\frac{(y-y^*)^2 + (z-z^*)^2}{2}\bigg]\\
\leq & \bigg[\frac{2sy^*}{a({x^*}^2 +a)} - 1\bigg] (x - x^*)^{2}  + \frac{(cx^* - b({x^*}^2 + a))}{a\beta d} (y - y^*)^2 \\
&+\frac{s}{a^2}\bigg[\frac{(x-x^*)^2 + (y-y^*)^2}{2}\bigg] \\
+& \bigg[\frac{q}{4br\alpha(q - p(\beta+\frac{\beta}{4b}+r))} - \frac{{x^*}^2 +a}{a\beta(\beta+\frac{\beta}{4b}+d)}\bigg]\bigg[\frac{(y-y^*)^2 + (z-z^*)^2}{2}\bigg]\\
=& \bigg[\frac{2sy^*}{a({x^*}^2 +a)} + \frac{s}{2a^2}- 1\bigg] (x - x^*)^{2} + \bigg[\frac{cx^* - b({x^*}^2 + a)}{a\beta d} + \frac{s}{2a^2}+ \\ &\frac{1}{2}\bigg(\frac{q}{4br\alpha(q - p(\beta+\frac{\beta}{4b}+r))} - \frac{{x^*}^2 +a}{a\beta(\beta+\frac{\beta}{4b}+d)}\bigg)\bigg] (y-y^*)^2 \\
&+ \frac{1}{2}\bigg[\frac{q}{4br\alpha(q - p(\beta+\frac{\beta}{4b}+r))} - \frac{{x^*}^2 +a}{a\beta(\beta+\frac{\beta}{4b}+d)}\bigg](z-z^*)^2.
\end{split}
\end{eqnarray}
One can note that $^c_{0}D^{m}_{t}V(x,y,z) \leq0, \forall (x,y,z)\in \Re^{3}_{+}$ if each coefficient of $(x-x^*)^2$, $(y-y^*)^2$ and $(z-z^*)^2$ are negative, giving the conditions
\begin{eqnarray}\nonumber
\begin{split}
&(i)~ \frac{2sy^*}{a({x^*}^2 +a)} + \frac{s}{2a^2}- 1 &< 0,\\
& (ii)~ \frac{cx^* - b({x^*}^2 + a)}{a\beta d} + \frac{s}{2a^2} \\
&+ \frac{1}{2}\bigg(\frac{q}{4br\alpha(q - p(\beta+\frac{\beta}{4b}+r))} - \frac{{x^*}^2 +a}{a\beta(\beta+\frac{\beta}{4b}+d)}\bigg) &< 0, \\
& (iii)~ \frac{q}{4br\alpha(q - p(\beta+\frac{\beta}{4b}+r))} - \frac{{x^*}^2 +a}{a\beta(\beta+\frac{\beta}{4b}+d)} &< 0.
\end{split}
\end{eqnarray}
Here $^c_{0}D^{m}_{t}V(x,y,z) = 0$ implies that $(x,y,z) = (x^*, y^*, z^*)$. Therefore, the only invariant set on which $^c_{0}D^{m}_{t}V(x,y,z) = 0$ is the singleton set $\{E^*\}$. Then, using Lemma (\ref{Global}) in \cite{HET15}, it follows that the interior equilibrium $E^*$ is global asymptotically stable for any $m \in (0,1]$. Hence the theorem is proven.\\

\noindent \textbf{Remark 1} This global stability result is independent of fractional order $m$ and it is also true for integer order $(m = 1)$.

\section{Numerical Simulations}
In this section, we perform extensive numerical computations of the fractional order system (\ref{Tritophic fractional order model}) for different fractional values of $m$ $(0 < m < 1)$ and also for $m=1$. We use Adams-type predictor corrector method (PECE) for the numerical solution of system (\ref{Tritophic fractional order model}). It is an effective method to give numerical solutions of both linear and nonlinear FODE \cite{Diethelm02,Diethelm04}. We first replace our system (\ref{Tritophic fractional order model}) by the following equivalent fractional integral equations:
\begin{eqnarray}\label{Tritophic fractional integral eqn} \nonumber
X(T) & = & X(0) + D^{-m}_{T} [a_0 X - b_0 X^2 - \frac{v_0 XY}{d_0 + X^2}], \nonumber\\
Y(T) & = & Y(0) + D^{-m}_{T} [-a_1 Y + \frac{v_1 XY}{d_1 + X^2} - \frac{v_2 YZ}{d_2 + Y}], \\
Z(T) & = & Z(0) + D^{-m}_{T} [c_3 Z^2 - \frac{v_3 Z^2}{d_3 + Y}]. \nonumber
\end{eqnarray}
and then apply the PECE (Predict, Evaluate, Correct, Evaluate) method.\\

\noindent Several examples are presented to illustrate the analytical results obtained in the previous section. Specially, our main objective is to explore the possibility of dynamical behavior of the fractional order system (\ref{Tritophic fractional order model}) by depending on the sensitive parameter and as well as the fractional order by keeping others parameters unchanged. To understand the effect of fractional order on the system dynamics, We varied $m$ in its range $0< m<1$. We also plotted the solutions for $m=1$, whenever necessary, to compare the solution of fractional order system with that of integer order. In numerical simulations, Initial values are indicated with stars and equilibrium points are denoted by red circles. \\

\textbf{Example 1.} In this example, here the parameter values are chosen as $b_0 = 0.075,~ a_1 = 0.105,~ d_1 = d_2 = 10.0,~ d_3 = 20.0,~ v_0 = 1.0,~ v_1 = 2.0,~ v_2 = 0.405,~ v_3 = 1.0$ and the initial condition $(1.2, 1.2, 1.2)$. All the parameters are taken from \cite{Ali16}. The bifurcation diagram with respect to sensitive parameters $a_0$ and $c_3$ is shown in Fig. 1 for different fractional order $m = 0.95, 0.75$ and the standard order $m = 1$. For the standard order $m = 1$, it is observed that the system (\ref{Tritophic fractional order model}) approches to chaos via period doubling bifurcation for $a_0 \in (0.25, 0.5)$ and $c_3 = 0.047$ (see Fig. 1(a)). It is interesting to note that the bifurcation disappears slowly with the decreasing of fractional order $m$ (see Figs. 1(b) and 1(c)).
\begin{figure}[H]
	\centering
	\includegraphics[width=10in, height=2.5in]{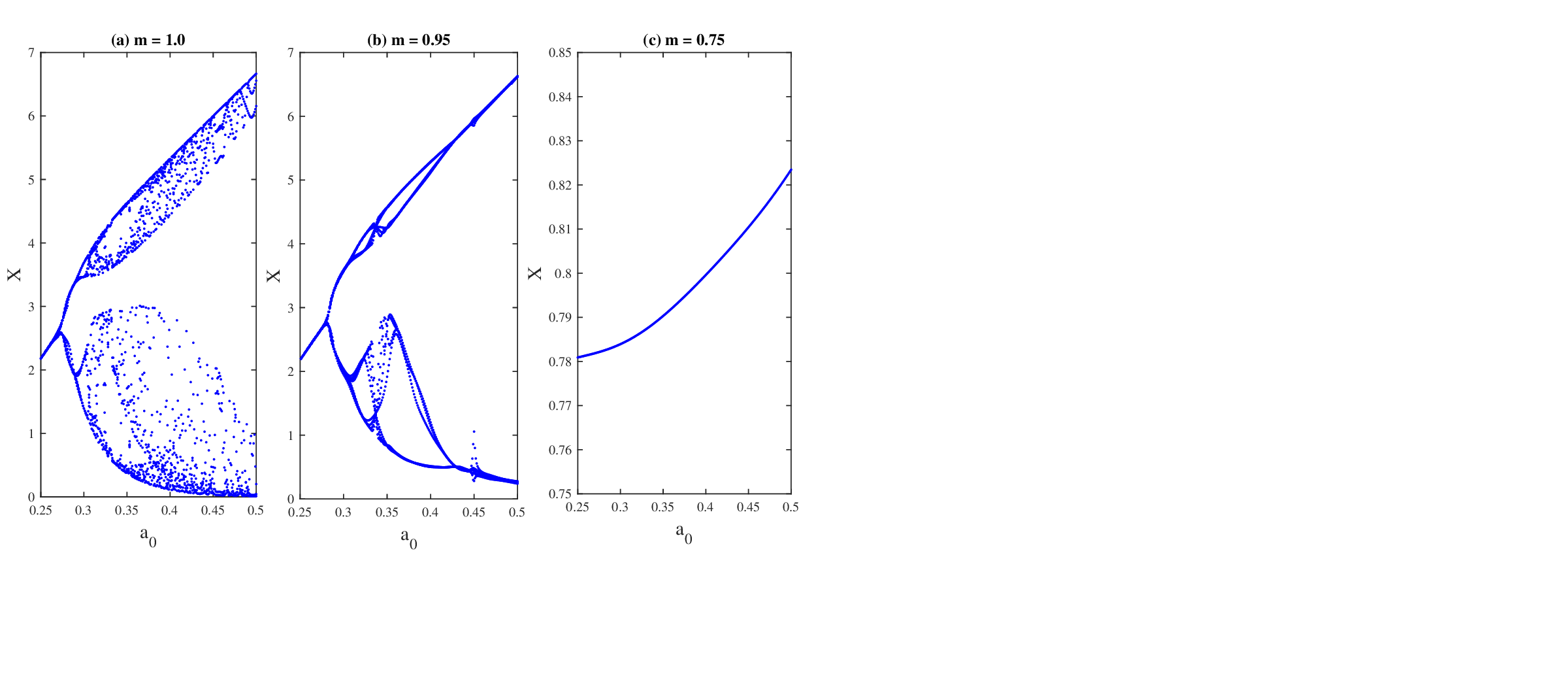}
	\vspace{-0.5in}
	\caption{Bifurcation diagram of system (\ref{Tritophic fractional order model}) for the $X$ population with respect to $a_0$ in $(0.25, 0.5)$ with different fractional orders $ m = 0.95,~ 0.75$ (Fig. 1(b) and 1(c)) and integer order $m = 1$ (Fig. 1(a)). Here $b_0 = 0.075,~ a_1 = 0.105,~ d_1 = d_2 = 10.0,~ d_3 = 20.0,~ v_0 = 1.0,~ v_1 = 2.0,~ v_2 = 0.405,~ v_3 = 1.0$ with $c_3 = 0.047$.}
\end{figure}
\textbf{Example 2.} Here we fixed $a_0 = 0.47$ (say) and varying $c_3 \in (0.041, 0.049)$ and keeping remaining parameters unaltered as in example 1. Both time series and phase portrait of our system (\ref{Tritophic fractional order model}) have been presented for different fractional order $m  = 0.95, 0.75$ and the standard order $m = 1$. In this case, we observe that chaotic behavior of our system changes to stability with decreasing of fractional order $m$ (see Figs 2(a) - 2(f)).

\begin{figure}[H]
	\centering
	\includegraphics[width=10in, height=2.5in]{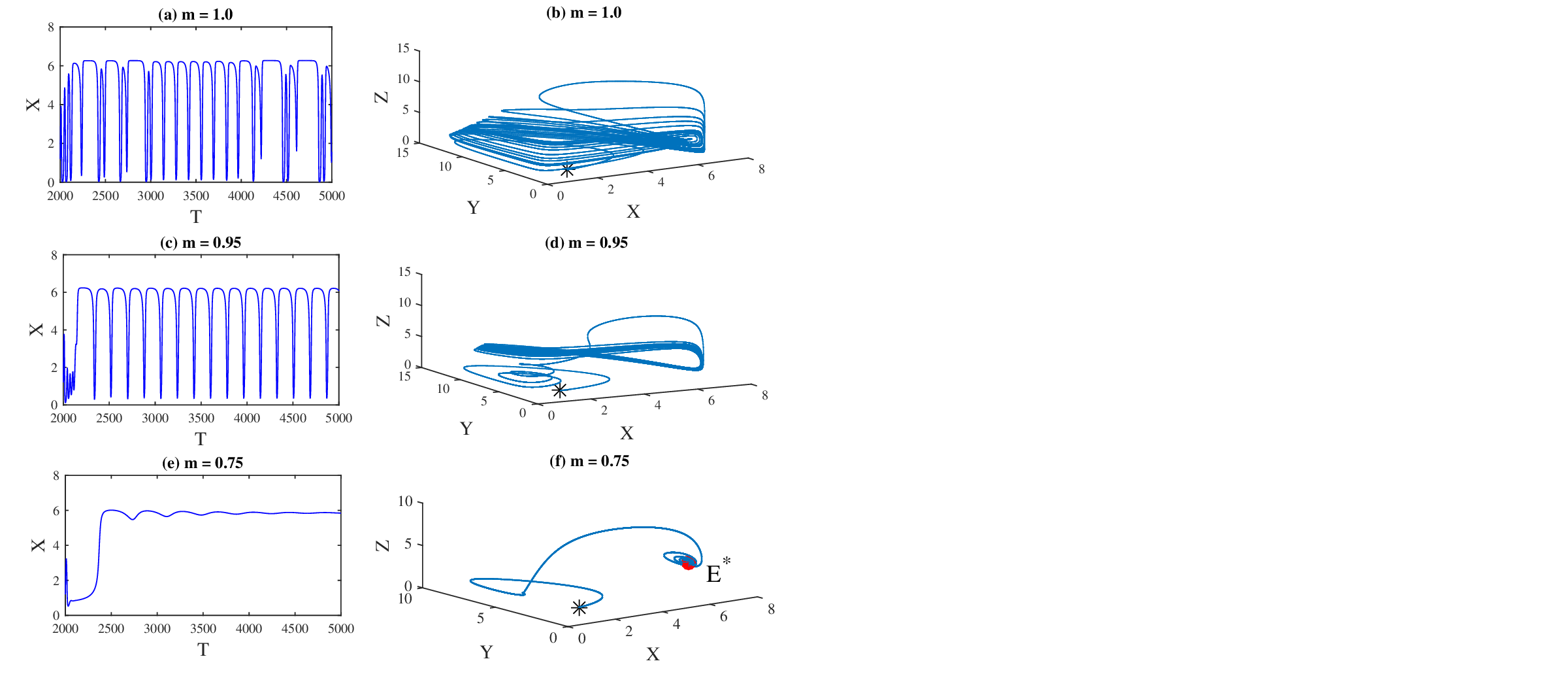}
	\vspace{-0.2in}
	\caption{The trajectory and phase portrait of system (\ref{Tritophic fractional order model}) with different fractional orders $ m = 0.95,~ 0.75$ (Fig. 2(c) - 2(f)) and integer order $m = 1$ (Fig. 2(a) - 2(b)). We observe that unstable behavior of our system changes to stability with decreasing of fractional order $m$. All the parameters are same as in example 1 with $a_0 = 0.47$ and $c_3 = 0.047$.}
\end{figure}

\textbf{Example 3.} Keeping $c_3$ unaltered, here we choose a smaller value of $a_0 = 0.27$ (say) and remaining all parameters are taken from example 1. Initial values are indicated with stars and equilibrium values are denoted by red circles in the figure. Step size for all simulations is considered as $0.05$. Using the above parameter set, we first verify the existence criteria of $E^*$. Here we observe $y^* - \frac{a}{s} = -1.4644<0, ~b - \frac{cx^*}{a +{x^*}^2} = -0.7582 < 0$ and $v_3 - c_3 d_3 = 0.06 >0$. Hence $E^* = (2.5772, 1.2766, 5.7002)$ exists in $R^3_+$ . Then compute $D(F) = -0.0084< 0$, $ A_{1} = 0.4033 >0$, $A_{2} =  0.0689 >0$, $A_{3} = 0.0221>0$. Thus, following Theorem (\ref{Interior Stability}) $(ii)$, the interior equilibrium $E^{*} = (x^*, y^*, z^*) = (0.7158, 1.3134, 8.7966)$ should stable for $0 < m < \frac{2}{3}$. In Fig. 3, we plot the time series solutions and phase portrait of FDE system (\ref{Tritophic fractional order model}) with different values of $m = 0.65,~ 0.60 < \frac{2}{3}$. It shows that all populations remain stable for all values of $m < \frac{2}{3}$, though solutions reach to equilibrium value more slowly as the value of $m$ becomes smaller (see Figs. 3(a) - 3(d) ).

\begin{figure}[H]
	\hspace{0.5in}
	\includegraphics[width=10in, height=3.5in]{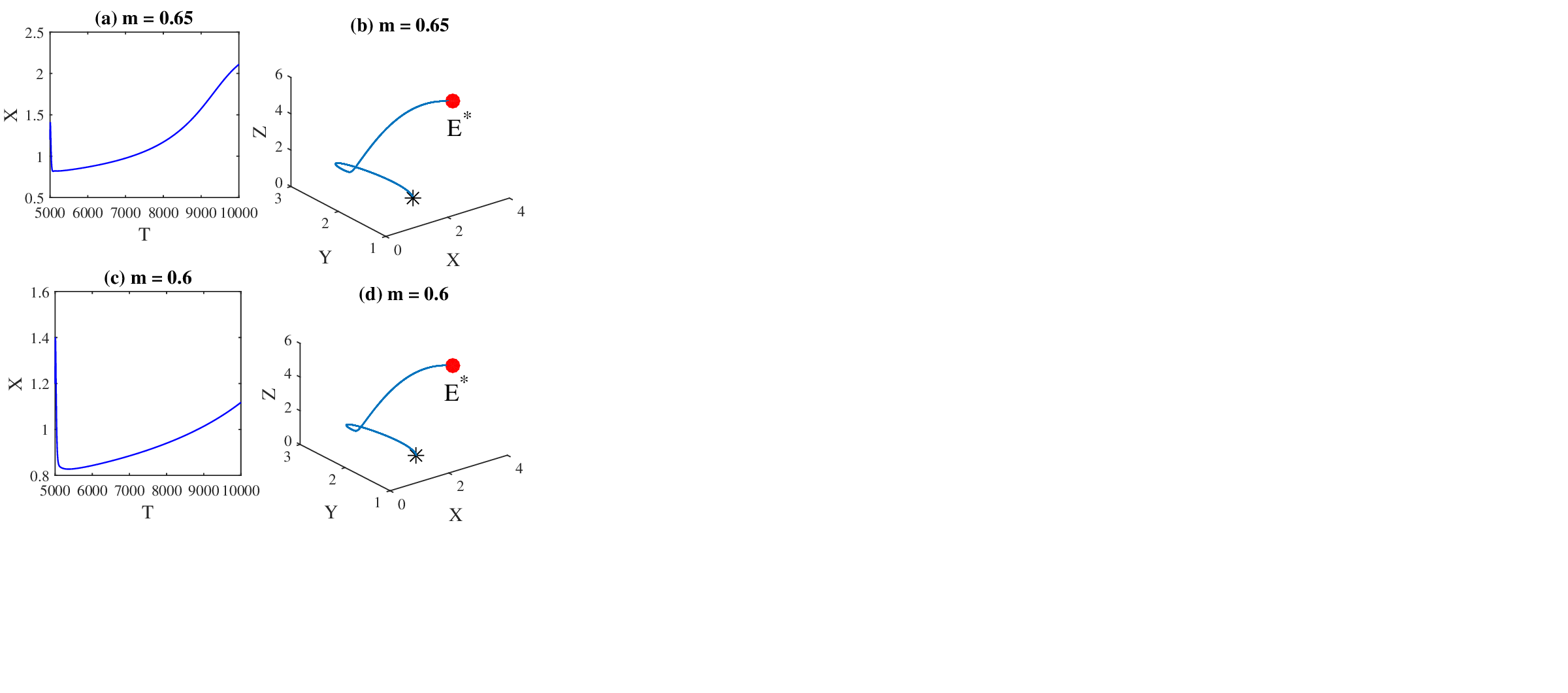}
	\vspace{-1in}
	\caption{The trajectory and phase portrait of system (\ref{Tritophic fractional order model}) with different fractional orders $ m = 0.65,~ 0.60 < \frac{2}{3}$ (Fig. 3(a) - 3(d)). We observe that the solution converges to interior equilibrium point for any values of $m < \frac{2}{3}$. It reaches to equilibrium value more slowly as the value of $m$ becomes smaller. All the parameters are same as in example 1 with $a_0 = 0.27$ and $c_3 = 0.047$.}
\end{figure}

Again if we increase the value of $a_0 = 0.35$ and keeping all parameters unaltered as in example 1, we see that our system (\ref{Tritophic fractional order model}) exhibits 2-periodic limit cycle, 1-periodic limit cycle for higher values of fractional order $m = 0.85$ as well as for integer order $m = 1$ (see Figs. 4(a) - 4(d)). If we decrease the value of $m$, then limit cycle disappears and system becomes stable. Here we choose $m = 0.75$ and observe that solution converges to interior equilibrium point $E^* = (4.0150, 1.2766, 0.7816, 5.6362)$ (see Figs. 4(e) - 4(f)).\\
\begin{figure}[H]
	\centering
	\includegraphics[width=10in, height=2.5in]{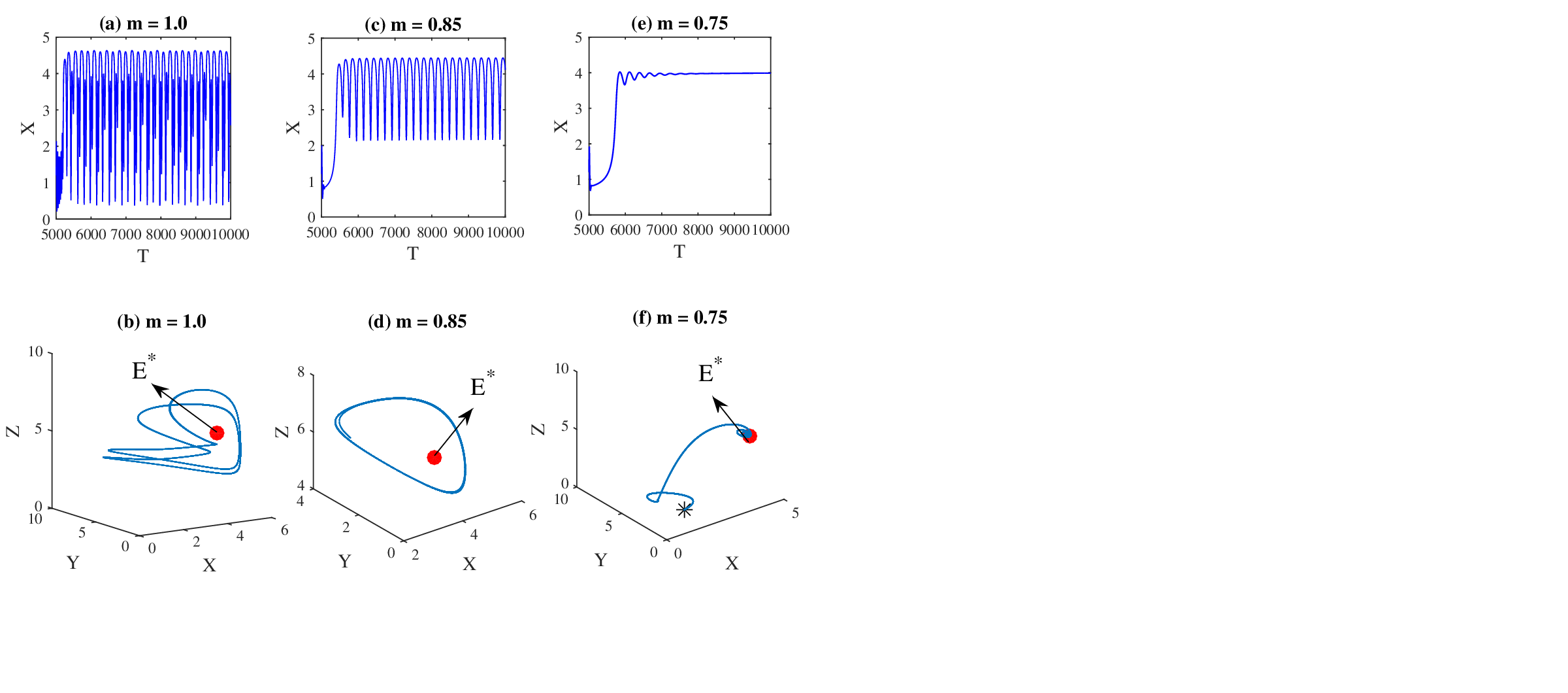}
	\vspace{-0.5in}
	\caption{The trajectory and phase portrait of system (\ref{Tritophic fractional order model}) with different fractional orders $ m = 0.85,~ 0.75$ (Fig. 4(c) - 4(f)). We observe that the solution converges to interior equilibrium point for any values of $m < \frac{2}{3}$. It reaches to equilibrium value more slowly as the value of $m$ becomes smaller. All the parameters are same as in example 1 with $a_0 = 0.35$ and $c_3 = 0.047$.}
\end{figure}

\textbf{Example 4.} Again we choose another parameter set $b_0 = 0.03,~ a_1 = 0.001,~ c_3 = 0.047,~ d_1 = d_2 = 10.0,~ d_3 = 20.0,~ v_0 = 0.85,~ v_1 = 2.5,~ v_2 = 2.2,~ v_3 = 1.0$ and  keeping same initial condition as in example 1, here we also choose a smaller value of $a_0 = 0.15$. Using the above parameter set, we first verify the existence criteria of $E^*$. Here we observe $y^* - \frac{a}{s} = -0.5532<0, ~b - \frac{cx^*}{a +{x^*}^2} = -2.6280 < 0$ and $v_3 - c_3 d_3 = 0.06 >0$. Hence $E^* = (3.2296, 1.2766, 2.0205)$ exists in $R^3_+$ . Then compute $D(F) = -0.0217< 0$, $ A_{1} = -0.0131<0$, $A_{2} =  -0.0044<0$. Thus, following Theorem (\ref{Interior Stability}) $(iii)$, the interior equilibrium $E^{*}$ should unstable for $m > \frac{2}{3}$. In Fig. 5, we plot the time series for $Y$ population and draw phase portrait of FDE system (\ref{Tritophic fractional order model}) on $XY$ plane with different values of $m = 0.95,~ 0.85 > \frac{2}{3}$. It shows that $Y$ population become unstable for different values of $m > \frac{2}{3}$ and our system (\ref{Tritophic fractional order model}) exhibits limit cycle around $E^*$. (see Figs. 5(a) - 5(d) ).
\begin{figure}[H]
	\centering
	\includegraphics[width=10in, height=2.5in]{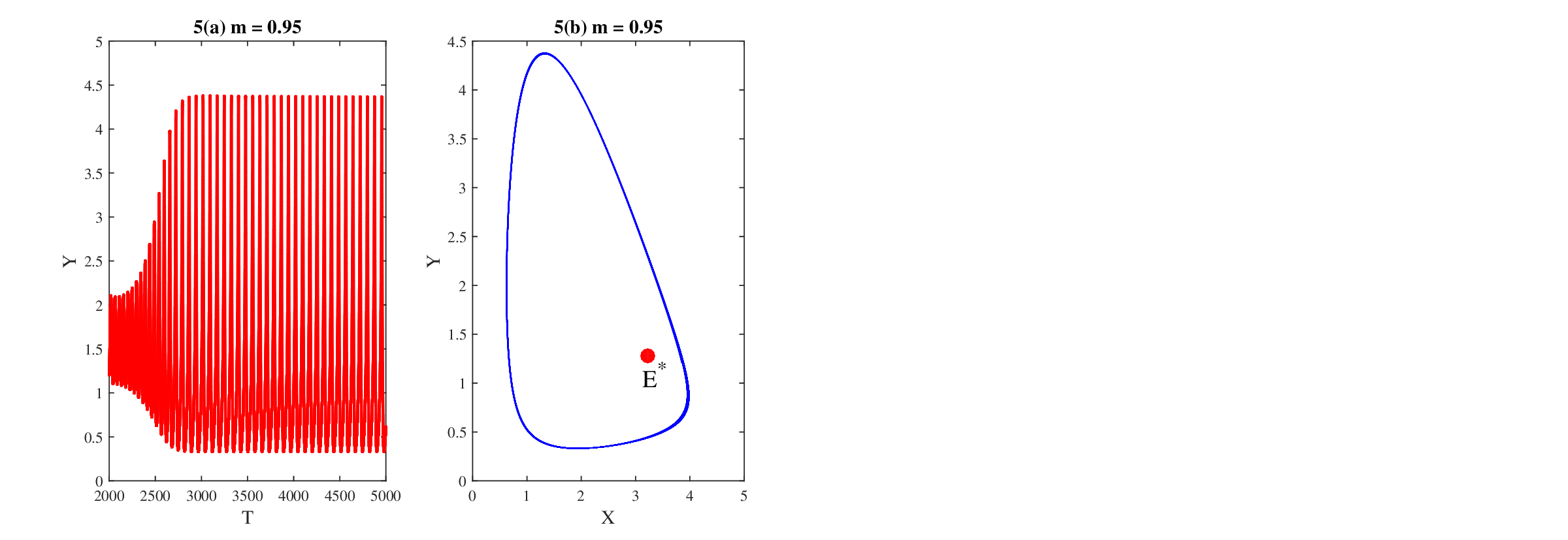}
\end{figure}
\begin{figure}[H]
	\centering
	\includegraphics[width=10in, height=2.5in]{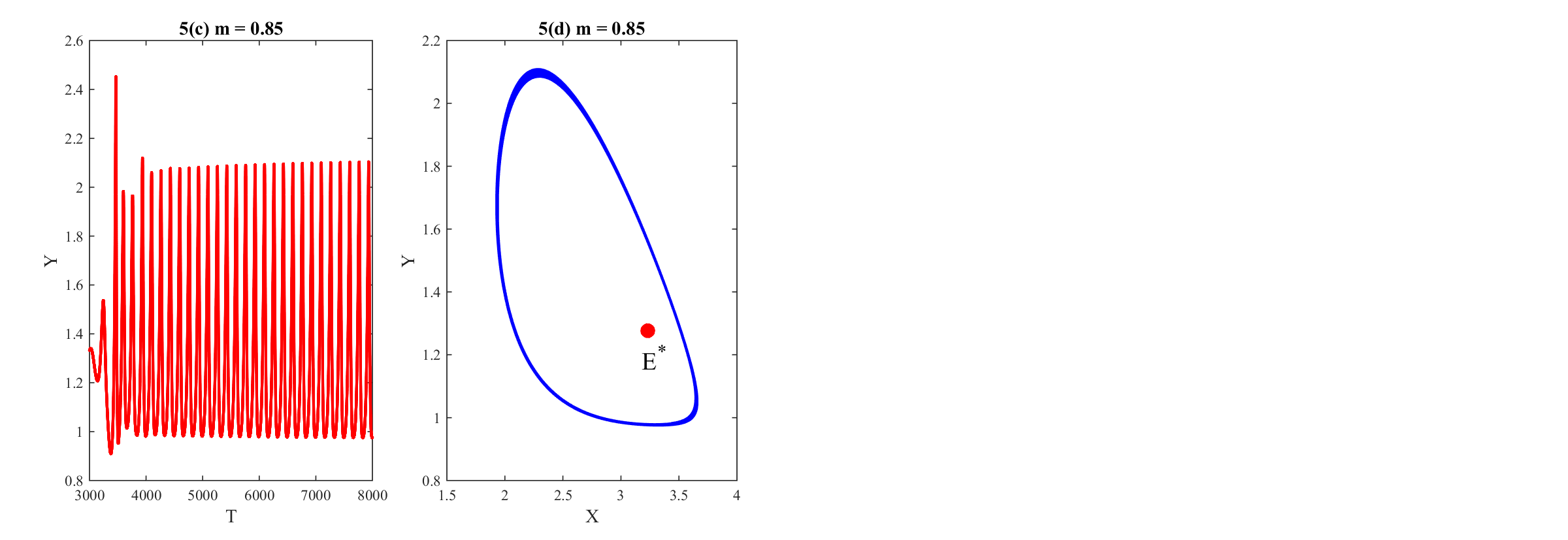}
	\caption{The trajectory and phase portrait of system (\ref{Tritophic fractional order model}) on $XY$ plane with different fractional orders $ m = 0.95,~ 0.85 > \frac{2}{3}$ (Fig. 5(a) - 5(d)). We observe that the $Y$ population becomes unstable for different values of $m > \frac{2}{3}$.  Here $b_0 = 0.03,~ a_1 = 0.001,~ c_3 = 0.047,~ d_1 = d_2 = 10.0,~ d_3 = 20.0,~ v_0 = 0.85,~ v_1 = 2.5,~ v_2 = 2.2,~ v_3 = 1.0$ and initial values are same as in example 1 with $a_0 = 0.15$.} 
\end{figure}

\noindent \textbf{Example 5:} To demonstrate the global stability of the interior equilibrium point $E^*$, we consider the parameter values $a_0 = 0.47$, $b_0 = 0.25$, $v_0 = 1.0$, $d_0 = d_1 = d_2 = 10.0$, $a_1 = 0.105$, $v_1 = 2.0$, $v_2 = 0.405$, $v_3 = 1.0$, $c_3 = 0.047$, $d_3 = 20.0$ and different initial points $(1.2, 1.2, 1.2)$, $(5.1,2.1,3)$, $(3,1,5)$, $(2,5,3.5)$, $(3,1,2)$, $(2.5,5,4)$, $(1.5,5.5,2)$, $(4.5,5.5,5)$. Initial values are indicated with stars and equilibrium values are denoted by red circles in the figure. Step size for all simulations is considered as $0.05$. Using the above parameter set, we first verify the existence criteria of $E^*$. Here we observe $y^* - \frac{a}{s} = -3.8744<0, ~b - \frac{cx^*}{a +{x^*}^2} = -0.2884 < 0$ and $v_3 - c_3 d_3 = 0.06 >0$. Hence $E^* = (1.4589, 1.2766, 3.7751)$ exists in $R^3_+$ . With these parameter values, we verify that all conditions of Theorem (\ref{Global_Stability}) are satisfied as 
\begin{eqnarray}\nonumber
\begin{split}
&(i)~ \frac{2sy^*}{a({x^*}^2 +a)} + \frac{s}{2a^2}- 1 & = -0.8084 < 0, \\
&(ii)~\frac{cx^* - b({x^*}^2 + a)}{a\beta d} + \frac{s}{2a^2} \\
& + \frac{1}{2}\bigg(\frac{q}{4br\alpha(q - p(\beta+\frac{\beta}{4b}+r))} - \frac{{x^*}^2 +a}{a\beta(\beta+\frac{\beta}{4b}+d)}\bigg) & = -0.0906 < 0,\\
&(iii)~\frac{q}{4br\alpha(q - p(\beta+\frac{\beta}{4b}+r))} - \frac{{x^*}^2 +a}{a\beta(\beta+\frac{\beta}{4b}+d)} &= -0.2623 < 0.
\end{split}
\end{eqnarray}
where $\alpha = \frac{1}{b^2 (c +\frac{c}{4b} +r)} = 0.6330>0$.  Fig. 5 demonstrates that solutions starting from different initial values converge to the equilibrium point $E^{*} = (1.4589, 1.2766, 3.7751)$ for different fractional orders, $m = 0.65, 0.75, 0.85$, and also for the integer order, $m = 1$, depicting the global stability of the interior equilibrium point for fractional order as well as integer order (see Figs. 6(a) - 6(d)).

\begin{figure}[H]
	\centering
	\includegraphics[width=10in, height=2.5in]{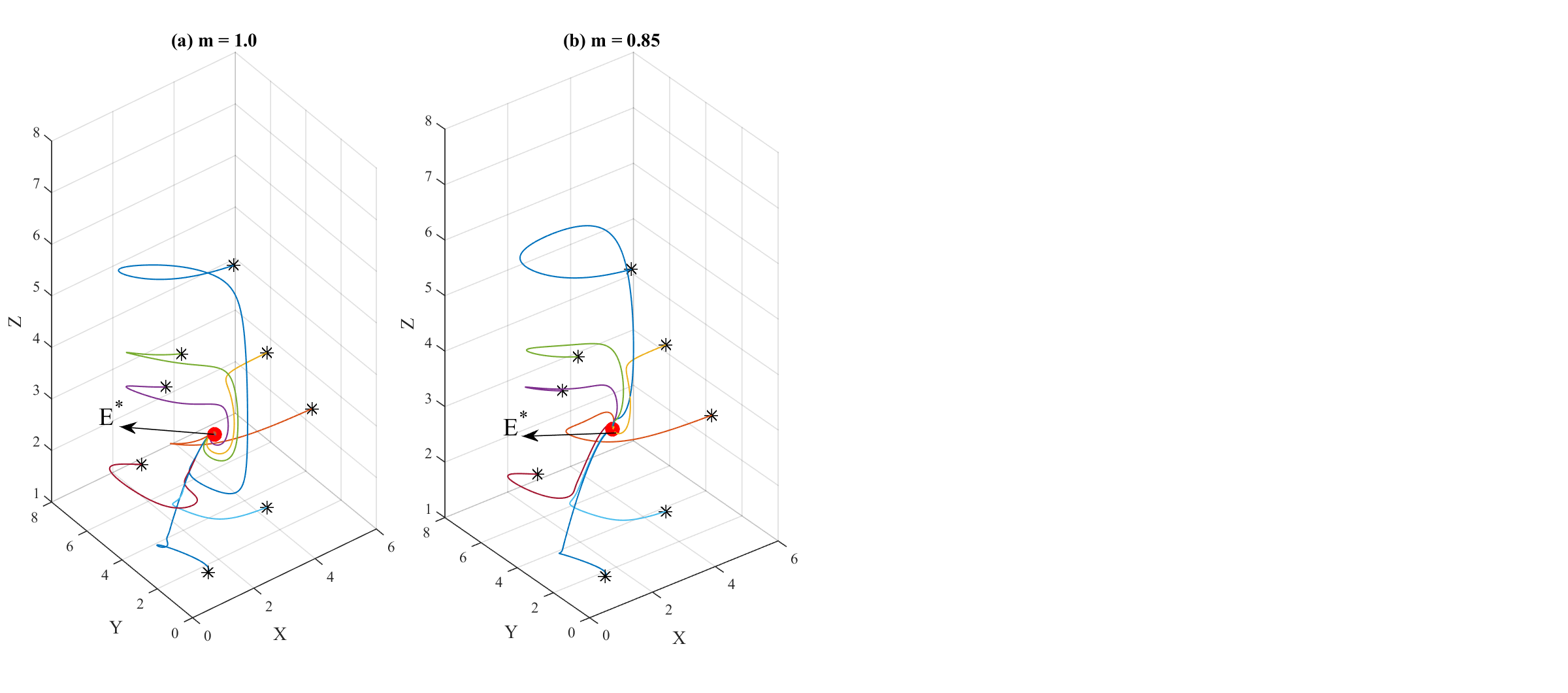}
\end{figure}
\begin{figure}[H]
	\centering
	\includegraphics[width=10in, height=2.5in]{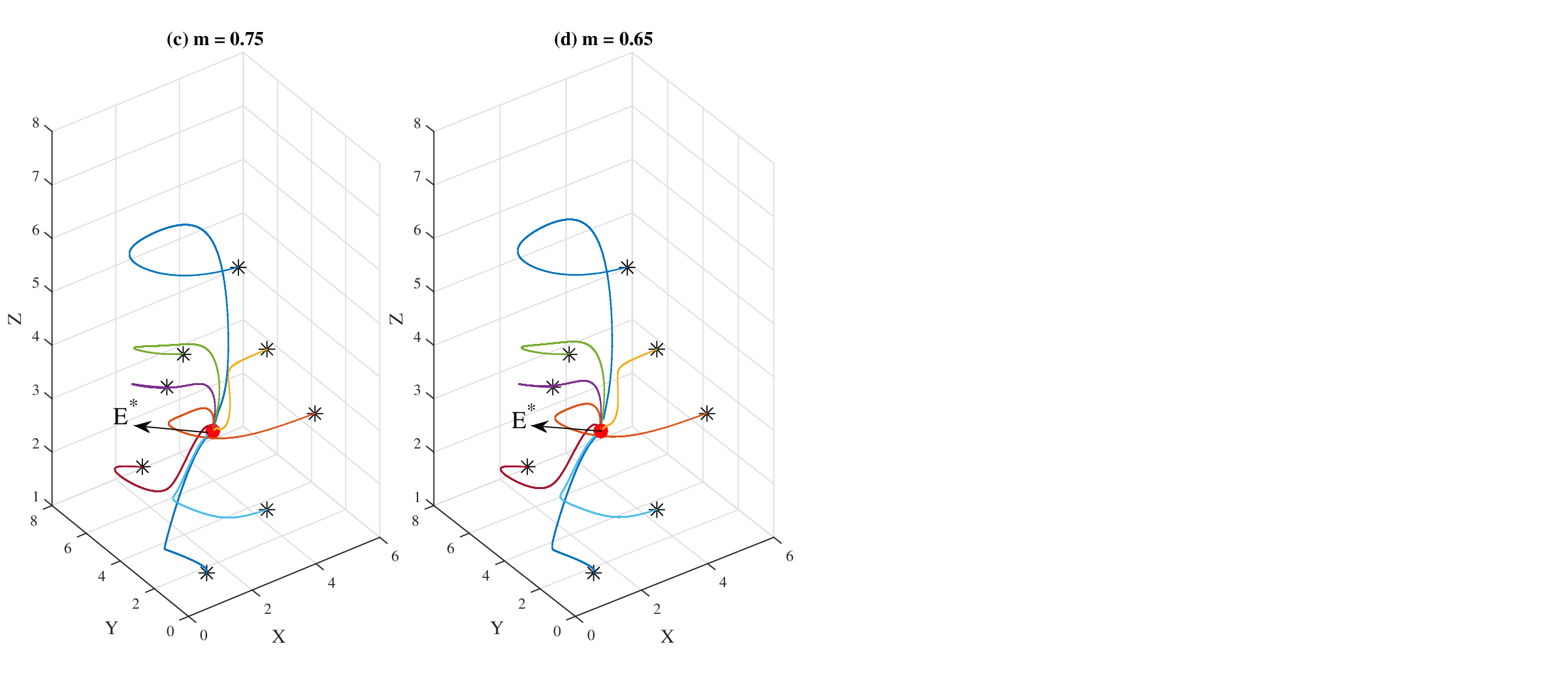}
	\caption{Trajectories with different initial values converge to the interior equilibrium point $E^*$ for different values of $m$, indicating global stability of the equilibrium $E^*$ when conditions of Theorem (\ref{Global_Stability}) are satisfied. All parameters are as in Fig. 1 except $b_0 = 0.25$.} 
\end{figure}

\section{Conclusions} \label{sec:6}

This paper generalizes the results of integer order three species food chain model \cite{Ali16} with simplified Holling type IV functional response. Here we first constructed a fractional order three species food chain model considering the fractional derivatives in Caputo sense. We proved different mathematical results like positivity and boundedness, existence and uniqueness,  local stability of different equilibrium points. It has been also proved that our system (\ref{Tritophic fractional order model}) is dissipative for any fractional order $0<m\leq 1$. Global stability of the interior equilibrium point have been only discussed. We defined suitable Lyapunov function to prove that the interior equilibrium is globally asymptotically stable if the system parameters satisfy some conditions. In such a case, the system does not show any complicated dynamics like chaos, indicating its global stability for any fractional order $0<m \leq 1$. This is more reinforced by the fact that solutions initiating from biologically feasible arbitrary initial points converge to the interior equilibrium point. To confirm the analytical results of our system, numerical simulation is performed for different sets of biologically feasible parameter values. Simulation results also agree perfectly with the analytical results. Numerically, it has been observed that the fractional order system shows more complex dynamics, like chaos as fractional order becomes larger and shows more simpler dynamics as the order $m$ decreases and becomes stable for lower value of $m$. Moreover, dynamics of the fractional-order system not only depends on system parameters but also depends on fractional order $m$.\\



\end{document}